\documentclass[11pt]{amsart}

\usepackage{amsmath}
\usepackage{amsthm}
\usepackage{amssymb}
\usepackage{amscd}
\usepackage{amsfonts}
\usepackage{amsbsy}
\usepackage{graphicx}

\usepackage{hyperref}
\hypersetup{colorlinks=true,citecolor=blue,%
        linkcolor=blue,urlcolor=blue,bookmarksnumbered=true,
        bookmarksopen=true,bookmarksopenlevel=1,breaklinks=true}

\def\epsilon{\varepsilon}
\newtheorem{theorem}{Theorem}[section]
\newtheorem{proposition}[theorem]{Proposition}
\newtheorem{lemma}[theorem]{Lemma}
\newtheorem{corollary}[theorem]{Corollary}

\numberwithin{equation}{section}

\numberwithin{figure}{section}

\begin{document}

\title[Centers, invariant straight lines and its configurations]
{Centers and invariant straight lines of planar real  polynomial vector fields and its configurations}

\author[H. He]
{Hongjin He}

\address{School of Mathematical Sciences, CMA-Shanghai, Shanghai Jiao Tong University, Shanghai 200240, China} \email{hehongjin000@126.com}

\author[C. Liu]
{Changjian Liu$^{\dag}$}
\address{School of Mathematics (Zhuhai), Sun Yat-sen University, 519086, Zhuhai, China}
\email{liuchangj@mail.sysu.edu.cn}

\author[D. Xiao]
{Dongmei Xiao$^{\ddag}$}
\address{School of Mathematical Sciences, CMA-Shanghai, Shanghai Jiao Tong University, Shanghai 200240, China} \email{xiaodm@sjtu.edu.cn}

%\thanks{$^{*}$  The work was partially supported by National Key R $\&$ D Program of China (No. 2022YFA1005900).}

%\thanks{($^{\dag}$) The author was supported by the NSFC grants (No. 12171491).}
%No.11371248).}

\thanks{($^{\ddag}$) Corresponding author.}

%\subjclass[2010]{34C07, 34C08, 32B10}



\maketitle
\begin{abstract}In the paper, we first give the least upper bound formula on the number of centers of planar real polynomial Hamiltonian vector fields.
  This formula reveals that the greater the number of invariant straight lines of the vector field and the less the number of its centers. Then we obtain some rules on the configurations
  of centers of planar real polynomial Hamiltonian Kolmogorov vector fields when the number of centers is exactly the least upper bound. As an application
  of these results,  we give an affirmative answer to a conjecture on the topological classification of configurations for the cubic Hamiltonian Kolmogorov vector fields with four centers. Moreover, we discuss the relationship between the number of centers of planar real polynomial vector fields and the existence of limit cycles, and prove that cubic real polynomial Kolmogorov vector fields have no limit cycles if the number of its centers reaches the maximum. More precisely, it is shown that the cubic real polynomial Kolmogorov vector field must have an elementary first integral in $\mathbb{R}^2\setminus\{xy=0\}$ if it has four centers, and the number of configurations of its centers  is one more than that of the cubic polynomial Hamiltonian Kolmogorov vector fields.
\end{abstract}

%\tableofcontents

\section{Introduction}

Hilbert 16th problem
 has two parts, see \cite{Hi}. The first part  is mainly to ask the relative position of the closed and  separate branches  ({\it ovals}, briefly) of real algebraic curves  of the $n$-th order
 in the real projective plane $\mathbb{RP}^2$ when their number is the maximum $(n-1)(n-2)/2+1$,  determined by Harnack. The second part is to ask the maximum number and relative position of limit cycles for
 planar real polynomial differential equations
\begin{equation}\label{generalfg}
 \frac{dx}{dt}=f(x,y),\ \ \frac{dy}{dt}=g(x,y), \ (x,y)\in \mathbb{R}^2,
 \end{equation}
 where $f(x,y)$ and $g(x,y)$ are real coefficients polynomial of degree $n$ and $m$ in two real variables $x$ and $y$, respectively.  And the limit cycle of  system \eqref{generalfg} is an isolated closed orbit (looks like oval) in $\mathbb{R}^2$.

It is well-known that the second part of the sixteenth problem
remains unsolved even for planar quadratic differential equations, see \cite{Il,IE,Li, Sm1, Sm2, Viro} and references therein.  An inherent complexity of this problem is implied by the fact that the limit cycles of polynomial  differential equations \eqref{generalfg} are not algebraic curves usually. Fortunately, the answer to the first part of Hilbert 16th problem can provide  some  hints for the study of the second part, for example, the ovals
in $\mathbb{R}^2$ of a $n$-th order real algebraic curve  can become limit cycles (called {\it algebraic limit cycles}) of some $n$-th order real polynomial differential equations \eqref{generalfg} in $\mathbb{R}^2$, and the number and  relative  position of the ovals in $\mathbb{R}^2$ is  that of the
algebraic limit cycles (see \cite{Chr} and references therein).
Thereby, one of the most fundamental questions about studying  Hilbert 16th problem is to understand the topological classification of real algebraic curves of fixed degree polynomials in $\mathbb{R}^2$.

Let $H(x,y)$ be a polynomial with real coefficients of degree $n+1$ in two real independent variables $x$ and $y$. Then
$$
X=-\frac{\partial H}{\partial y}\frac{\partial }{\partial x}+\frac{\partial H}{\partial x}\frac{\partial }{\partial y}=(-\frac{\partial H}{\partial y}, \frac{\partial H}{\partial x})
$$
is a planar polynomial Hamiltonian vector field of degree $n$. A critical point $p$ of $H(x,y)$ is called {\it center} of $X$ if there exists a neighborhood filled with periodic orbits (ovals) of $X$ with the exception of the critical point $p$. This neighborhood is usually called {\it period annulus} of $X$.  {\it Hilbert's 16th Problem on a period annulus} is to study how many  limit cycles can bifurcate from the families of ovals  by a small polynomial perturbation, see \cite{Ar, Ber, Bin, ChLi, FGX, FHX, Ga,Ili, HoIli} and references therein. And the relative position of the limit cycles of the perturbed vector fields is closely related to the number and configuration of centers of $X$. So the number of centers of $X$ and its configurations play an important role in study the second part of Hilbert's 16th problem.  However, the complete configurations of the centers of $X$ have been known only for $n=1, 2$, see \cite{Ga}, even thought there have been many interesting results for several subclasses of cubic polynomials vector fields, see for instance \cite{CLV1, CLV2, CLV3, LX, RS, VS, Zo2, Zo3} and references therein.

This paper has two purposes. One is to study the number of centers of $X$ and its configurations when their number is the maximum, which  is the first step for the study of the number and the relative position of limit cycles generated by perturbating Hamiltonian systems with centers. And the other is to study dynamics of planar real polynomial vector fields with the maximum centers, which reveals some relationship between the number of centers  and the existence of limit cycles for the vector fields.

To our knowledge,  Cima, Gasull and Ma$\tilde{n}$osas first studied the maximum number of centers,  gave a beautiful  upper bound of the number of centers of planar polynomial Hamiltonian vector fields $X$ depending on its degree in \cite{CGM} and  the upper bound  can be achieved by some $X$ with no  critical point at infinity or only a pair of critical points at infinity, where the critical point at infinity can be defined by using the Poincar\'e compactification of vector fields. A  critical point at infinity is called an {\it infinite critical point} of the vector field. A question arises naturally  how many is the least upper bound on the number of centers for any a polynomial Hamiltonian vector field $X$ with $r$ pairs of infinite critical points? where $r$ is a nonnegative integer. Moreover,  how many is the
possible configurations of centers if the number of its centers is the upper bound?  Motivated by these questions, we study the maximum number of centers of any planar polynomial Hamiltonian vector fields $X$ and  possible configurations of centers.

  The first goal of this paper is to answer the question and find the  least upper bound formula of the number $c(X)$ of centers of $X$, which  depends on the degree of $X$ and the number of  its infinite critical points. This improves and generalizes
 the main theorem in \cite{CGM}. In particular, the existence of an invariant straight line implies that a pair of infinite critical points of $X$ exists.  Accordingly, this least upper bound formula establishes the connection the number of centers of planar polynomial Hamiltonian vector fields  with the number of its invariant straight lines. This will help us to obtain configurations of the centers of $X$.
  In \cite{LX} authors studied  configurations of centers for cubic polynomial Hamiltonian vector fields  with two intersecting invariant straight lines,  and conjectured that there are only two types of configurations of centers if this cubic vector field has four centers. In
  the paper,  we consider general polynomial Hamiltonian vector fields with only two intersecting invariant straight lines, that is,  polynomials $H(x,y)$ with real coefficients of degree $n+1$ have only two different linear factors.  Without loss of generality, we assume that the polynomials can  be written as $H(x,y)=xyF(x,y)$, where  $F(x,y)$ are real polynomials of degree $n-1$ with $n-1\ge 2$. Hence, the corresponding Hamiltonian vector fields are
$$
 X_{hk}=\left (-x(F(x,y)+y\frac{\partial F}{\partial y}), y(F(x,y)+x\frac{\partial F}{\partial x})\right ),
$$
 which was called {\it polynomial Hamiltonian Kolmogorov vector fields} (or {\it HK-vector fields} for short) in \cite{LX}.

 The second aim in the paper is to study the possible configurations of centers of  $X_{hk}$ if the number of its centers  is exactly the  least upper bound. We obtain some rules on configurations of centers of $X_{hk}$ by index theory and perturbation technique. Especially, when $n=3$, using these rules we can prove the  conjecture in \cite{LX} is true, that is, the cubic HK-vector fields  with four centers  have only two different types of configurations of centers. Moreover,  we also describe completely the different possible global phase portraits of this vector field in Poincar\'e disk.

The last purpose of this paper is to discuss  whether planar real polynomial  system \eqref{generalfg} has no limit cycles if the number of its centers is the maximum. The problem on quadratic polynomial vector fields has been solved, see \cite{cLi,Sc,Vu,Zo,Zo1} and references therein. However, this problem remains unsolved for cubic polynomial vector field.
And  the least upper bound on the number of centers is still open for the  vector fields of associated system \eqref{generalfg}
$$Y=f(x,y)\frac{\partial }{\partial x}+g(x,y)\frac{\partial }{\partial y}=(f(x,y), g(x,y)),$$
 where $f(x,y)$ and $g(x,y)$ are real coefficients polynomials of degree $n$ and $m$, respectively if ${\rm max}\{m, n\}\ge 4$, see \cite{Gas}. Note that a critical point $p$ of $Y$ is  a real solution $(x_0,y_0)$ of $f(x,y)=0$ and $g(x,y)=0$, which is also called {\it finite critical point of $Y$}.  We would like to connect the existence of a first integral  in $\mathbb{R}^2\setminus \mathcal{B}$ and the maximum number of its centers for system \eqref{generalfg}, where the set $\mathcal{B}$ consists of these orbits whose limit sets contains only  critical points of system \eqref{generalfg}.  Of course, system \eqref{generalfg} may have a first integral defined in $\mathbb{R}^2\setminus \mathcal{B}$ even thought the number of its centers is not the maximum, for example, system \eqref{generalfg} has a global center in \cite{HLX,HX}.

 In this paper, we study the problem for planar cubic polynomial vector fields with two intersecting invariant straight lines.  Without loss of generality,  the cubic polynomial vector fields with two intersecting invariant straight lines can be written as
 $$
  Y_{k}=(xP(x,y), yQ(x,y)),
 $$
 where $P(x,y)$  and $Q(x,y)$ are all quadratic real polynomials. We call $Y_k$ as {\it cubic polynomial Kolmogorov vector fields}. It is clear that $Y_k$ has at most four centers.
Combining  different techniques from algebraic curves, topology and  differential equations, we  prove that $Y_k$  has no limit cycles if the number of its centers reaches the maximum. More precisely, it is shown that the cubic real polynomial Kolmogorov vector field must have an elementary first integral in $\mathbb{R}^2\setminus\{(x,y):\ xy=0\}$ if it has four centers, and there are only three different kinds of configurations of centers of $Y_k$ in some equivalent sense.
 This reveals that non-Hamiltonian cubic polynomial Kolmogorov vector fields with four centers can have one more configurations of centers than that of Hamiltonian ones.

\smallskip

In the study, we mainly use some properties of algebraic curves in the complex projective plane $\mathbb{CP}^2$,  Poincar\'e compactification of planar polynomial vector fields and index theory, and develop some perturbation techniques and qualitative analytical method.
This paper is organized as follows. In sake of convenience for readers, in Section 2 we introduce some notations and necessary preliminaries from algebraic curves and vector fields, and provide some known results in literature. In Section 3, we obtain the  least upper bound formula of the number $c(X)$ of centers for planar polynomial Hamiltonian vector fields $X$.  In Section 4, we investigate the possible configurations of centers of  HK-vector fields $X_{hk}$ if the number of its centers is exactly the least upper bound. As an application,  in last section  we study  dynamics of the cubic Kolmogorov vector fields $X_{hk}$ and $Y_{k}$ when they have four centers, respectively. The conjecture proposed in \cite{LX} is proved, and the difference on the configurations of the four centers is discovered for cubic polynomial Hamiltonian  and  non-Hamiltonian Kolmogorov vector fields.

\section{Preliminaries }\label{notation}

In the section, we first introduce some notations and concepts in algebraic curves and vector fields, then review some known results (for detail please see the literature \cite{BK, CGM, CGM2, CL, Ful} and references therein). All of them  will be used later.

%\subsection{The concepts  from algebraic curves and vector fields}

Let us denote  the set of all real planar polynomial vector fields $(f(x,y), g(x,y))$  by $\mathcal{Y}_{n,m}$, without loss of generality, we always assume $n\geq m$ in this article. Then
  $$
  \mathcal{Y}_{n,m}=\left \{Y:Y=(f(x,y), g(x,y)),\; {\rm deg}(f)=n, {\rm deg}(g)=m, n\ge m,\ (x,y)\in \mathbb{R}^2 \right\}.
  $$
Notice that $f(x,y)$ and $g(x,y)$ can be expanded into the sum of homogeneous polynomials as follows.
\begin{equation}\label{fj}
f(x,y)=\sum_{i=0}^nf_i(x,y),\; \; g(x,y)=\sum_{j=0}^mg_j(x,y),
\end{equation}
where $f_i(x,y)$ and $g_j(x,y)$ are the $i$th order and $j$th order homogeneous parts of $f(x,y)$ and $g(x,y)$, respectively, $i=0,1,\cdots,n$ and $j=0,1,\cdots,m$.

Using the Poincar\'e compactification of $Y$, we can calculate the infinite critical points of $Y$, and obtain that  $(x_0,y_0)$  is {\it an infinite critical point of $Y$} if and only if it is {\it a nonzero real solution} of the following $n$th order homogeneous polynomial equation
\begin{equation}\label{infiniteCRq}
-yf_n(x,y)+ xg_n(x,y)=0 \ {\rm if}\ n=m,
\end{equation}
or
\begin{equation}\label{infiniteCRfq}
  -yf_n(x,y)=0 \ {\rm if}\ n>m.
\end{equation}
Clearly, the infinite critical points of $Y$ appear in pairs of
diametrally opposite points, see \cite{G}.

%In the following we recall some concepts in algebraic curves.
From the algebraic viewpoint,  we usually  discuss the common zeros of polynomials $f(x,y)$ and $g(x,y)$ in the complex plane $\mathbb{C}^2$ and the complex projective plane $\mathbb{CP}^2$. Any a point in $\mathbb{CP}^2$ can be represented by using its projective coordinates in three local charts $W_i$, $i=1,2,3$,
$$\mathbb{CP}^2=\cup_{i=1}^3 W_i,\ W_i=\{[x:y:z]:=[x_1:x_2:x_3]\in\mathbb{CP}^2,\  x_i\neq 0\}.$$%\cong \mathbb{C}^2.$$
Thus, we identify each $(x,y)\in \mathbb{C}^2$ with $[x:y:1]\in \mathbb{CP}^2$, and points
at infinity of $\mathbb{C}^2$ on the straight line are of the form $[x:y:0]\in \mathbb{CP}^2$ with $x^2+y^2\not=0$. We often use  coordinates $(x,y,z)$ to instead  homogeneous coordinates $[x:y:z]$ if there is no confusion. Therefore,
the common zero $p=(x_*,y_*,0)$ of polynomials $f(x,y)$ and $g(x,y)$ is at infinity of $\mathbb{C}^2$ if and only if
\begin{equation}\label{infiniteCP}
f_n(x_*,y_*)=0=g_m(x_*,y_*).
\end{equation}
Obviously, $(x_*,y_*)$ is a nonzero solution of \eqref{infiniteCRq} or \eqref{infiniteCRfq} if $(x_*,y_*)\in \mathbb{R}^2$. Hence, a common zero $p=(x_*,y_*,0)$ of polynomials $f(x,y)$ and $g(x,y)$  at infinity of  $\mathbb{C}^2$ is an infinite critical point of the vector field $Y$ if $(x_*,y_*)\in \mathbb{R}^2$. But an infinite critical point of the vector field $Y$ may not be a  common zero  of polynomials $f(x,y)$ and $g(x,y)$  at infinity of  $\mathbb{C}^2$ by \eqref{infiniteCRq} or \eqref{infiniteCRfq}.

Note that two real polynomials $f(x,y)$ and $g(x,y)$  have no common components in $\mathbb{R}[x,y]$ if and only if $f(x,y)$ and $g(x,y)$  have no common components in $\mathbb{C}[x,y]$.
Hereafter we consider the vector field $Y\in\mathcal{Y}_{n,m}$ in the complex number domain $\mathbb{C}$ or complex plane $\mathbb{C}^2$ for convenience. And the original real vector fields $Y\in\mathcal{Y}_{n,m}$ can be regarded as this complex vector fields confined in $\mathbb{R}^2$.

Let $f(x,y)\in\mathbb{C}[x,y]$ be a complex coefficients polynomial of degree $n$ in two variables $x$ and $y$. The {\it affine plane curve} $f(x,y)$ is the zero set of this polynomial  $$V(f):=\{(x,y)\in\mathbb{C}^2:\ f(x,y)=0\}.$$
Let $p=(x_0,y_0)$ be a point in $V(f)$. A natural number $k$ is called {\it the multiplicity of the curve $f(x,y)$ at $p$}, denoted by {\it $ k=m_p(f)$}, if
$$f(x+x_0,y+y_0)=f_k(x,y)+f_{k+1}(x,y)+\cdots +f_n(x,y), \ f_k(x,y)\not\equiv 0, $$
where $f_i(x,y)$ is the $i$th order homogeneous polynomial in two variables $x$ and $y$, $k\le i\le n$. Clearly $k\ge 1$. Then there exist natural numbers $r_i\in \mathbb{N}$ and complex numbers $a_i,b_i\in \mathbb{C}$ such that
$$f_k(x,y)=\prod_{\sum r_i=k} L_i^{r_i},\ L_i=a_ix+b_iy, \ L_i\not=L_j \ \rm{for}\ i\not=j. $$
The line $L_i$ in $\mathbb{C}^2$ is called {\it the tangent lines to $f(x,y)$ at $p$} and $r_i$ is called {\it the multiplicity of this tangent $L_i$}.

Let $\mathcal{O}_p(\mathbb{C}^2)$ be the ring of rational functions defined at a point $p\in \mathbb{C}^2$. And let $<f,g>\mathcal{O}_p(\mathbb{C}^2)$ be the ideal generated by two affine plane curves $f(x,y)$ and $g(x,y)$ in $\mathcal{O}_p(\mathbb{C}^2)$. Then {\it the intersection number of $f(x,y)$ and $g(x,y)$} at the point $p$, denoted by $I(p,f\cap g)$,  is defined by
$$
I(p,f\cap g)=dim_\mathbb{C}\frac{\mathcal{O}_p(\mathbb{C}^2)}{<f,g>\mathcal{O}_p(\mathbb{C}^2)}.
$$
The intersection number is the unique number which satisfies some properties (see \cite{Ful} for detail), and these properties also tell us how to calculate the intersection number of $f(x,y)$ and $g(x,y)$ at the point $p$.

Consider homogenization $f^*$ of $n$-th order polynomial $f\in\mathbb{C}[x,y]$,
$$f^*=z^nf(\frac{x}{z},\frac{y}{z})=\sum_{i=0}^{n}z^{n-i}f_i(x,y).$$
$f$ is regard as the restriction of the projective plane algebraic curve $f^*$ on the chart
$$W_3=\{(x,y,1):\ (x,y,z)\in W_3, z=1\}$$ since $f^*(x,y,1)=f(x,y)$. Therefore, for a point $p=[x_0:y_0:1]\in \mathbb{CP}^2$ and two projective plane algebraic curves $f^*$ and  $g^*$,  it can be defined {\it the intersection number of  $f^*$ and  $g^*$ at $p$},  $I([x_0:y_0:1], f^*\cap g^*)$,  by $I((x_0,y_0),f\cap g)$. It is not hard to prove that the intersection number $I([x_0:y_0:1], f^*\cap g^*)$ does not depend on the choice of charts. A well-known result about the intersection number of two projective plane algebraic curves is B\'ezout's theorem as follows, whose proof can be found in \cite{BK, Ful}.

{\bf B\'ezout Theorem}
{\it Let $f$ and $g$ be projective plane algebraic curves of degree $n$ and $m$ respectively. If $f$ and $g$ do not have common components, and $A$ is the set of all common zeros of $f$ and $g$ in $\mathbb{CP}^2$. Then
$$
\sum_{p\in \mathbb{CP}^2}I(p,f\cap g)=\sum_{p\in A}I(p,f\cap g)=nm.
$$}

We now consider such a vector field $Y\in\mathcal{Y}_{n,m}$ whose elements $f(x,y)$ and $g(x,y)$ {\it  have no common zeros at infinity} in  $\mathbb{CP}^2$, which is equivalent to the condition that $f_n(x,y)$ and $g_m(x,y)$ do not have common components in $\mathbb{C}^2$. Let us denote the subset consisting of these vector fields by $\Psi_{n,m}$, that is,
   $$
   \Psi_{n,m}=\left\{Y: Y\in \mathcal{Y}_{n,m},  f_n(x,y)\ {\rm and}\ g_m(x,y)\ {\rm  have \ no\ common \ components} \right\}.
   $$
It can be check that $\mathcal{Y}_{n,m}\backslash\Psi_{n,m}$ is contained in an algebraic hypersurface of $\mathcal{Y}_{n,m}$, that is,  {\it $\Psi_{n,m}$ is generic in $\mathcal{Y}_{n,m}$}, see \cite{CL2} for detail.

Further we consider vector fields $Y\in \mathcal{Y}_{n,m}$ such that $f(x,y)$ and $g(x,y)$ {\it have exactly $nm$ different common zeros} in  $\mathbb{C}^2$, denoted the set consisting of these  vector fields by $G_{n,m}$,
$$
G_{n,m}=\{Y:\ Y\in \mathcal{Y}_{n,m} \ {\rm and}\ \sharp\{(x,y)\in\mathbb{C}^2:f(x,y)=g(x,y)=0\}=nm\}\},
$$
where $\sharp\{\cdot\}$ denotes the number of elements of the set $\{\cdot\}$.

Then by B\'ezout's theorem,  $G_{n,m}\subset \Psi_{n,m}$. And these common zeros $p=[x_0:y_0:1]$ of $f(x,y)$ and $g(x,y)$ in $\mathbb{CP}^2$ are all {\it finite critical points} of the vector fields $Y\in G_{n,m}$ if  $(x_0,y_0)\in \mathbb{R}^2$. Thus, they are isolated and elementary, where the critical point $p$ of vector fields  $Y=(f(x,y), g(x,y))$ is called {\it elementary} if  Jacobian matrix of the vector field $(f,g)$ with respect to $(x,y)$ at $(x_0,y_0)$ has no zero eigenvalues. It is easily proved that  $\mathcal{Y}_{n,m}\backslash G_{n,m}$ is also contained in an algebraic hypersurface of $\mathcal{Y}_{n,m}$, which implies that {\it $G_{n,m}$ is generic in $\mathcal{Y}_{n,m}$} too.

Note that vector field $Y$ in $\mathbb{R}^2$ can be induced two vector fields $\tilde{Y}_{\pm}$ in the northern hemispheres $\mathbb{S}_+^2$ and southern hemispheres $\mathbb{S}_-^2$, respectively via central projection by considering the plane $\mathbb{R}^2$ as the tangent space at the north pole of the unit sphere $\mathbb{S}^2$, called {\it  of Poincar\'e sphere}, in $\mathbb{R}^3$. The induced vector field $\tilde{Y}$ in each hemisphere is analytically conjugate to $Y$ in $\mathbb{R}^2$, and the equator $\mathbb{S}^1$ of $\mathbb{S}^2$ is bijective correspondence with the points at infinity of $\mathbb{R}^2$. The global dynamics of $Y$ in the whole $\mathbb{R}^2$ including its dynamical behavior near infinity is analytically conjugate to that of $\tilde{Y}$ in $\mathbb{S}_+^2\cup \mathbb{S}^1$, which is called {\it Poincar\'e disc}.  By using a scaling the independent time variable,  we can extend the induced vector fields $\tilde{Y}$ in $\mathbb{S}^2\setminus \mathbb{S}^1$ to an analytical vector field  $P(Y)$ defined in the whole  $\mathbb{S}^2$. This is the {\it Poincar\'e compactification} of $Y$ in $\mathbb{R}^2$, where $P(Y)$ is analytically equivalent to $Y$, and the analytical expression of $P(Y)$ can be computed in the six local charts $U_i$ of the differentiable manifold $\mathbb{S}^2$, see \cite{DLA} for detail. Then the  Poincar\'e-Hopf theorem tell us that the sum of the indices at the critical points of $P(Y)$ is equal to the Euler-Poincar\'e characteristic of the compact manifold $\mathbb{S}^2$ if  $P(Y)$ has only isolated critical points. Note that all  critical points of $P(Y)$ on $\mathbb{S}^2$ are isolated  if all  critical points of $\tilde{Y}$ in $\mathbb{S}_+^2\cup \mathbb{S}^1$  are isolated. And the indices of the corresponding critical points of $P(Y)$ and ${Y}$ are the same.

We use the same notations in \cite{CGM} to denote the sum of indices of all  isolated finite critical points (resp. all  isolated infinite critical points) of $Y$ by $\sum_f i$ (resp. $\sum_{inf} i$). Similarly, we can define the sum of the absolute values of indices of all isolated finite critical points (resp. all isolated infinite critical points) by $\sum_f |i|$ and $\sum_{inf} |i|$.  Hence, if $Y$ has finitely many critical points (including finite critical points and infinite critical points), then by Poincar\'e-Hopf theorem we have
\begin{equation}\label{phT}
2\sum_{f} i+\sum_{inf} i=2.
\end{equation}

 Now let us recall the index of a critical point $p$ of vector field $Y$ in $\mathbb{R}^2$. Assume that $\gamma$ is an oriented simple closed curve which does not pass through critical points of $Y$, and there is a unique critical point $p$ of $Y$ in the interior surrounded by $\gamma$. Then the topological degree of the map $h: \  \gamma \to S^1$ (the unit circle), given by $h(M)=\frac{Y(M)}{\|Y(M)\|}$ for $\forall M\in \gamma$,  is called {\it the index of the critical point $p$ of $Y$}, denoted by $i_Y(p)$. The index $i_Y(p)$ is an integer, which can be calculated by Poincar\'e method as follows:  given
 a direction vector $v$ in $\mathbb{R}^2$, we check if there exist only finitely many points $M_i\in \gamma$, $i=1,\cdots, k$ such that the direction of vector field $Y$ at the point $M_i$, denoted by $Y(M_i)$, is parallel to $v$.  Let $q_+$ (resp. $q_-$) be the number of points $M_i$ at which the vector $Y(M)$ passes through the given direction $v$ in the counterclockwise (resp. clockwise) sense when a point $M$ on $\gamma$ moves along  $\gamma$ in counterclockwise sense. Then the index $i_Y(p)$ of $p$ is
 $$ i_Y(p)=\frac{q_+ - q_-}{2},$$
see \cite{CL} for detail.

%\subsection{The estimation on the index of critical points of vector fields}

There have been some useful estimations on the index $i_Y(p)$. Let us revisit some of them in \cite{CGM, CL} which are used in this study.
\begin{lemma}\label{lem1}(Lemma 1.1 in  \cite{CL})
Let $p$ be an isolated critical point of a vector field $Y=(f(x,y),g(x,y))$. Then
$$
|i_Y(p)|\leq\min\{m_p f,m_p g\},
$$
where $m_p f$ and $m_p g$ are the multiplicity of algebraic curves $f(x,y)$ and $g(x,y)$ at $p$, respectively.
\end{lemma}
Note that the intersection number $I(p,f\cap g)$ has a property: $I(p,f\cap g)\geq m_p(f)m_p(g)$ and the equality holds if and only if $f(x,y)$ and $g(x,y)$ have no common tangent lines at $p$. By Lemma \ref{lem1} we have
\begin{lemma}\label{lem2}(Lemma 1.3 in \cite{CGM})
Let $p$ be an isolated critical point of a vector field $Y=(f(x,y),g(x,y))$. Then
$$
(i_Y(p))^2\leq I(p,f\cap g).
$$
\end{lemma}
Note that $I(p,f\cap g)\geq 1$ if $p\in f\cap g$. By B\'ezout's theorem, we can achieve an important estimation about the sum of the absolute values of indices of all isolated finite critical points.

\begin{proposition}\label{thm-est1}(Lemma 1.4 in \cite{CGM})
Assume that  a vector field  $Y\in \mathcal{Y}_{n,m}$. If  all  finite critical points of $Y$ are isolated, then
$$
\sum_f|i|\leq nm,
$$
\end{proposition}
Next result is the other important estimation of indices of all finite critical points, see appendix of \cite{CL} to get the proof.
\begin{proposition}\label{thm-est2}
Assume that  a vector field  $Y\in \mathcal{Y}_{n,m}$. If  all  finite critical points and all infinite critical points of $Y$ are isolated, then
$$
|\sum_f i|\leq \min\{n,m\}=m,
$$
\end{proposition}

Last we recall the relationship between the local dynamics  of Hamiltonian vector field $X$ at an isolated finite critical point and its index, and an estimation on maximum number of centers of polynomial vector field $Y\in \mathcal{Y}_{n,n}$ as follows.
\begin{lemma}\label{lem-f}(Proposition 2.1 in \cite{CGM})
Let $p$ be an isolated finite critical point of Hamiltonian vector field $X=(-\frac{\partial H}{\partial y},\frac{\partial H}{\partial x})\in\mathcal{Y}_{n,m}$. Then the index $i_X(p)\leq 1$ of $X$ at  $p$ characterizes the topology behaviour of orbits near $p$, i.e.
\begin{itemize}
\item[(i)] $i_X(p)=1$ if and only if the critical point $p$ is a center.
\item[(ii)] $i_X(p)=1-h\leq 0$ if and only if the neighbourhood of critical point $p$ is only composed of $2h$ hyperbolic sectors, where $h$ is a positive integer, and hyperbolic sector is saddle sector.
\end{itemize}
\end{lemma}

\begin{lemma}\label{lem-A}(Theorem A in \cite{CGM2})
Assume that $C_m$ is the maximum number of centers of polynomial vector field $Y\in \mathcal{Y}_{n,n}$, where $n>1$. Then
$$
[\frac{n^2+1}{2}]\le C_m\le \frac{n(n+1)}{2}-1,
$$
where $[\cdot]$ denotes the integer part of the number.
\end{lemma}

\section{The number of centers of Hamiltonian polynomial vector fields $X$ }\label{mnc}

In this section we consider Hamiltonian vector fields $X$ with polynomial Hamiltonian functions $H(x,y)$ of degree $n+1$, and polynomial $H_{n+1}(x,y)$ is the $(n+1)$-th order homogeneous parts of $H(x,y)$. Assume that $X\in \mathcal{Y}_{n,m}$, that is,
$${\rm deg}\left(-{\partial H(x,y)}/{\partial y}\right)=n,\ {\rm deg}\left({\partial H(x,y)}/{\partial x}\right)=m,\ n\ge m.$$ Then  $H_{n+1}(x,y)=ay^{n+1}$ with $a\neq 0$ if $n>m$.  From
\eqref{infiniteCRq} or \eqref{infiniteCRfq}, we know that the linear factors of $H_{n+1}(x,y)$ determine all infinite critical points of $X$. Hence, $X$ has finitely many infinite  critical points.
Assume that $H_{n+1}(x,y)$ has $r$ linear factors. Then $X$ has $r$ pairs of infinite  critical points, where $r$ is a nonnegative integer.
Our main result in the section is to give the least upper bound on number $c(X)$ of centers for polynomial Hamiltonian vector fields $X$ with exactly $2r$ infinite critical points as follows, which improves and generalizes Theorem 3.1 in \cite{CGM}.
\begin{theorem}\label{thmmain1}
Let $X=(-\frac{\partial H}{\partial y},\frac{\partial H}{\partial x})\in\mathcal{Y}_{n,m}$ be a polynomial Hamiltonian vector field, and let $c(X)$ be the number of centers of $X$. If vector fields $X$ have exactly $2r$ infinite critical points, then
$$
c(X)\leq C_{n,m}=\begin{cases}
&[\frac{n^2+1-r}{2}], \quad n=m,\\
&[\frac{nm+1}{2}], \quad n>m,
\end{cases}
$$
where $r$ is a nonnegative integer.
Moreover, this bound $C_{n,m}$ can be realized.
\end{theorem}

\subsection{The upper bound for $X$ with common components }

Before proving Theorem \ref{thmmain1}, we first prove an auxiliary result.
Note that $H(x,y)$ is any a polynomial  functions of degree $n+1$. So polynomials ${\partial H}/{\partial x}$ and $-{\partial H}/{\partial y}$  may have common components.
The following proposition provides an estimation of $c(X)$ if $X$ has non-isolated critical points.
\begin{proposition}\label{pro-non-iso}
If polynomials ${\partial H}/{\partial x}$ and $-{\partial H}/{\partial y}$ have common components, then the number $c(X)$ of centers of the polynomial Hamiltonian vector field $X$ satisfies that $c(X)\le C_{n,m}-1$.
\end{proposition}
\begin{proof}
Suppose that  polynomials ${\partial H}/{\partial x}$ and $-{\partial H}/{\partial y}$ have common factors $\bar{H}(x,y)$ which is a polynomial of degree $s$ ($1\leq s\leq m$). Then the corresponding Hamitionian system of $X$ can be written to
\begin{equation}\label{hcomf}
\begin{split}
\frac{dx}{dt}=&-\frac{\partial H}{\partial y}=-\bar{H}(x,y)\bar{f}(x,y),\\
\frac{dy}{dt}=&\frac{\partial H}{\partial x}=\bar{H}(x,y)\bar{g}(x,y),
\end{split}
\end{equation}
where $\bar{f}(x,y)$ and $\bar{g}(x,y)$ are polynomials of degree $n-s$ and $m-s$, respectively. And polynomials $\bar{f}(x,y)$ and $\bar{g}(x,y)$ have no common components in $\mathbb{R}^2$.

Consider the set
$$
B=\{(x,y): \ \bar{H}(x,y)=0, (x,y)\in \mathbb{R}^2\}\subset \mathbb{R}^2.
$$
Then $B$ has only two possibility:  $B=\emptyset$ or  $B\not=\emptyset$.
We now  study the number of centers of system \eqref{hcomf} in the two cases.

Case (i): $B=\emptyset$. Then either $\bar{H}>0$  or $\bar{H}<0$ in $\mathbb{R}^2$. Hence, system \eqref{hcomf} is orbitally equivalent to the following polynomial system
\begin{equation}\label{fgint}
\begin{split}
\frac{dx}{dt}=&-\bar{f}(x,y)=-\frac{1}{\bar{H}(x,y)}\frac{\partial H}{\partial y},\\
\frac{dy}{dt}=&\bar{g}(x,y)=\frac{1}{\bar{H}(x,y)}\frac{\partial H}{\partial x}
\end{split}
\end{equation}
by time scaling. Hence, system \eqref{hcomf} and system \eqref{fgint} have the same number of centers, that is $c(X)=c(\bar{X})$, here $\bar{X}=(-\bar{f}(x,y),\bar{g}(x,y))$.

Case (ii): $B\neq\emptyset$. Then system \eqref{hcomf} is  orbitally equivalent to  system \eqref{fgint} in each connected components of $\mathbb{R}^2\setminus B$ by time scaling. Hence, system \eqref{hcomf} and system \eqref{fgint} have the same number of centers in $\mathbb{R}^2\setminus B$.  Now we consider the critical points  of system \eqref{hcomf} in the set $B$. Suppose $p_0\in B$ is a center of system \eqref{hcomf}, we claim that $p_0$ is a center of system \eqref{fgint} too.

Indeed, if $p_0$ is a center of system \eqref{hcomf}, then $p_0$ must be an isolated zero of $\bar{H}(x,y)=0$ in $\mathbb{R}^2$. Thus, there exists a small neighbourhood $U(p_0)$ of $p_0$ such that either $\bar{H}(x,y)>0$ or $\bar{H}(x,y)<0$ in $U(p_0)\setminus \{p_0\}$. Hence, system \eqref{hcomf} is orbitally equivalent to system \eqref{fgint} in $U(p_0)\setminus\{p_0\}$ by time scaling, which implies that every orbits of system \eqref{fgint} are closed orbits in $U(p_0)\setminus\{p_0\}$.  By definition of center,  $p_0$ must be a center of system \eqref{fgint}. It follows that all  centers of system \eqref{hcomf} are the centers of system \eqref{fgint}. So
$$
c(X)\leq c({\bar{X}}).
$$

We now estimate the upper bound of $c({\bar{X}})$.  Since polynomials $\bar{f}(x,y)$ and $\bar{g}(x,y)$ do not have common components in $\mathbb{R}^2$, all finite critical points of system \eqref{fgint} are isolated. And the infinite critical points of system \eqref{fgint} correspond to the real linear factors of the polynomial
$$
x\bar{g}_{n-s}-y\bar{f}_{n-s}=\frac{1}{\bar{H}(x,y)}\left(x\frac{\partial H_{n+1}}{\partial x}+y\frac{\partial H_{n+1}}{\partial y}\right)=(n+1)\frac{H_{n+1}}{\bar{H}(x,y)}, \ {\rm as}\ n=m
$$
or
$$
-y\bar{f}_{n-s}=\frac{1}{\bar{H}(x,y)}\left(y\frac{\partial H_{n+1}}{\partial y}\right)=(n+1)\frac{H_{n+1}}{\bar{H}(x,y)}, \ {\rm as}\ n>m,
$$
where $\bar{f}_{n-s}$ and $\bar{g}_{n-s}$ are the highest homogenous parts of polynomials $\bar{f}(x,y)$ and $\bar{g}(x,y)$ respectively.
Thus, the infinite critical points of system \eqref{fgint} are isolated.
 By Lemma \ref{lem-f} and Proposition \ref{thm-est2}, we have
$$
c({\bar{X}})-\sum_{i_{\bar{X}}(p)\leq 0} i_{\bar{X}}(p)=\sum_f i\leq m-s.
$$
On the other hand, by Proposition \ref{thm-est1},
$$
c({\bar{X}})+\sum_{i_{\bar{X}}(p)\leq 0} i_{\bar{X}}(p)=\sum_f |i|\leq (n-s)(m-s).
$$
Hence, we have
$$
c({\bar{X}})\leq \frac{1}{2}(m-s+(n-s)(m-s))\le\frac{1}{2}n(m-1).
$$
If $n>m\geq 1$, then
$$c({\bar{X}})\leq\frac{1}{2}n(m-1)\leq\frac{1}{2}nm-1\leq[\frac{nm+1}{2}]-1.$$
If $n=m$, by Lemma \ref{lem-A},  we can obtain a better estimation
$$
c({\bar{X}})\leq \frac{1}{2}(n-s+(n-s)(n-s))-1\leq \frac{1}{2}n(n-1)-1.
$$
Note that $\frac{1}{2}n(n-1)\leq [\frac{n^2+1-r}{2}]$ for all $0\leq r\leq n+1$. Hence, $$c(X)\leq c({\bar{X}})\leq C_{n,m}-1$$ for all $n\geq m\geq 1$.
\end{proof}

\subsection{The least upper bound for $X$ with no common components}

From Proposition \ref{pro-non-iso}, we prove Theorem \ref{thmmain1} only in the case that $\frac{\partial H}{\partial x}$ and $\frac{\partial H}{\partial y}$ do not have common components. And we divide the proof of Theorem \ref{thmmain1} into two parts. The first part is to prove that the number $c(X)$ of centers has a upper bound $C_{n,m}$. And in the second part we prove the upper bound $C_{n,m}$ is sharp, that is, there exists a Hamiltonian vector field $X\in\mathcal{Y}_{n,m}$, which has $C_{n,m}$ centers.

\begin{proof}[Proof of the first part of theorem \ref{thmmain1}]
%In fact, we can suppose all the finite critical points of $X$ are isolated, otherwise it cannot get the maximum number of centers.

Since the degree of polynomials $\frac{\partial H}{\partial y}$ and $\frac{\partial H}{\partial x}$ is $n$ and $m$ respectively,  we distinguish two cases $n=m$ and $n>m$ to prove that  the number $c(X)$ of centers of  polynomial Hamiltonian vector fields $X$ has the upper bound $C_{n,m}$.

\textbf{Case 1.} when $n=m$, we consider Poincar\'e compactification $P(X)$ of $X$. We study the index of infinite critical points of $X$ in two local charts $U_1$ and $U_2$
$$U_1=\{(x,y,z)\in \mathbb{S}^2:\ x>0\},\quad U_2=\{(x,y,z)\in \mathbb{S}^2:\ y>0\}$$
of Poincar\'{e} sphere $\mathbb{S}^2$.
 Since the arguments on studying index are similar in the two local charts, without loss of generality, we assume that all infinite critical points of $X$ lie on the local chart  $U_1$. Thus,  those $r$ pairs of infinite critical points of $X$ are $p_i=(1,\pm y_i,0),\ i=1,2,\dots,r$. We claim
 that the index $i_{P(X)}(p_i)$ at infinite critical points $p_i$ of vector field $P(X)$ satisfies
\begin{equation}\label{eqi}
i_{P(X)}(p_i)\geq 1-I_i,
\end{equation}
where $I_i$ is the intersection number of polynomials $-\frac{\partial H}{\partial y}$ and $\frac{\partial H}{\partial x}$ at $p_i$ in $\mathbb{CP}^2$, that is,
$$
I_i=I(p_i,f^*\cap g^*)=I\left([1:y_i:0],\left(-\frac{\partial H}{\partial y}\right)^*\cap\left(\frac{\partial H}{\partial x}\right)^*\right).
$$
In fact, on the local chart $U_1$,  the corresponding differential system of $P(X)$ has the form
\begin{equation}\label{chartU1}
\begin{split}
\frac{dy}{dt}&=-yf^*(1,y,z)+g^*(1,y,z)=\sum_{i=0}^n(n+1-i)z^iH_{n+1-i}(1,y),\\
\frac{dz}{dt}&=-zf^*(1,y,z)=z\sum_{i=0}^nz^{i}\frac{\partial H_{n+1-i}}{\partial y}(1,y)
\end{split}
\end{equation}
by Poincar\'e transformation $x\mapsto \frac{1}{z},y\mapsto \frac{y}{z}$.
System  \eqref{chartU1} has critical points  $p_i=(\pm y_i,0), \ i=1, \cdots, r$ in $U_1$. To estimate the index of $P(X)$ at $p_i=(y_i,0)$,  by using Poincar\'e method, we choose a circle $S_\epsilon=\{(y,z):(y-y_i)^2+z^2=\epsilon^2,0<\epsilon\ll 1\}$  in $U_1$ such that $S_\epsilon$ does not pass through any critical points of \eqref{chartU1}. Given a direction vector $v=(1,0)$, it can be checked that the points on $S_\epsilon$ at which the direction vector of  \eqref{chartU1} being parallel to $v$ must satisfy $zf^*(1,y,z)=0$. We then consider a point $M$ on $S_\epsilon$ moving counterclockwise to calculate the number $q_+$ (resp. $q_-$) of points on $S_\epsilon\cap\{(y,z):zf^*(1,y,z)=0\}$ at which the direction vector of  \eqref{chartU1} at $M$ passes through the given direction $v$ in the counterclockwise (resp. clockwise) sense, that is, the number of times the direction vector $P(X)(M)$ passes through the given direction $v$ along $S_\epsilon$ in the counterclockwise (resp. clockwise) sense as follows.

 Firstly we calculate the contribution of points $S_\epsilon\cap\{z=0\}=\{(y_i-\epsilon,0),(y_i+\epsilon,0)\}$ to the index at $p_i$. Since $p_i$ is an infinite critical point of \eqref{chartU1}, it follows $H_{n+1}(1,y_i)=0$. Hence, $H_{n+1}(1,y)$ has the expression
$$
H_{n+1}(1,y)=a_l(y-y_i)^l+a_{l+1}(y-y_i)^{l+1} \dots+a_{n+1}(y-y_i)^{n+1},
$$
where $l\geq 1$ and $a_l\neq 0$. The direction vector of $P(X)$ at the point near $(y_i-\epsilon,0)$ and $(y_i+\epsilon,0)$ has the following approximation
\begin{equation*}
\begin{split}
P(X)&=((n+1)H_{n+1}(1,y)+O(z),z(\frac{\partial H_{n+1}}{\partial y}(1,y)+O(z)))\\
&\approx((n+1)a_l(y-y_i)^l,la_l(y-y_i)^{l-1}z).
\end{split}
\end{equation*}
Hence, this direction vector depends on the sign of $a_l$ and the parity of $l$.  After making some analysis on $a_l$ and $l$, we can sketch the four cases on direction vectors $P(X)$ at the point near $(y_i-\epsilon,0)$ and $(y_i+\epsilon,0)$ in counterclockwise sense for all $i=1,2,\dots,r$, see Figure \ref{fig1}.
\begin{figure}
  \centering

  \includegraphics[width=1\textwidth]{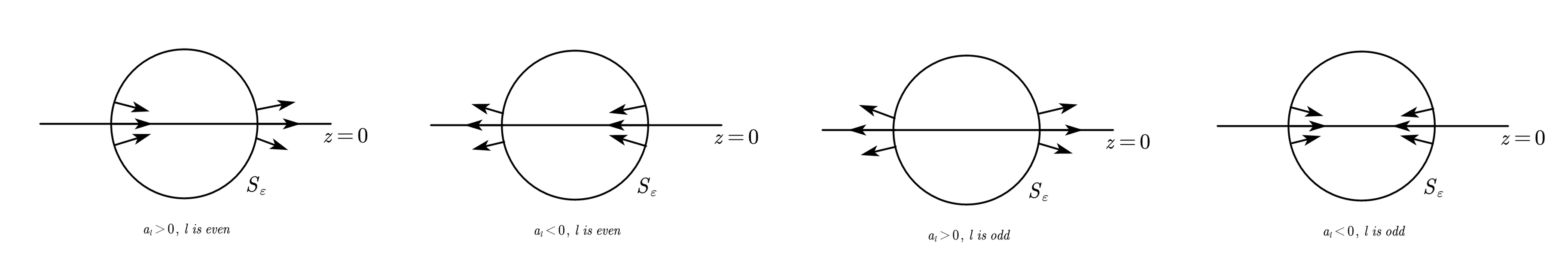}
  \caption{\small The direction vector $P(X)(M)$ on $S_\epsilon$ near the points $(y_i-\epsilon,0)$ and $(y_i+\epsilon,0)$.}
  \label{fig1}
\end{figure}

Next we consider the contribution of points $S_\epsilon\cap\{f^*(1,y,z)=0\}$ to the index at $p_i$.  Let ${q}_+^*$ (resp. ${q}_-^*$) be the number of points on curve $f^*(1,y,z)=0$ at which the direction vector $P(X)(M)$ pass through the direction $v$ along $S_\epsilon$ in the counterclockwise(resp. clockwise) sense. So we have
$$
q_+=2+{q}_+^*,\quad q_-={q}_-^*.
$$
Following the definition of index by Poincar\'e,  the index at point $p_i$ is
$$
i_{P(X)}(p_i)=\frac{q_+-q_-}{2}=1+\frac{{q}_+^*-{q}_-^*}{2}.
$$
Note that  the curve $f^*(1,y,z)=0$ and $S_\epsilon$ have at most $2m_{p_i}f^*$ common points when $0<\epsilon\ll 1$ (see \cite{Ber} for detail). So
$$
{q}_+^*+{q}_-^*\leq 2m_{p_i}f^*,
$$
where $m_{p_i}f^*$ is the multiplicity of the curve $f^*(1,y,z)$ at $p_i$.
Therefore,
\begin{equation}\label{leq1}
|i_{P(X)}(p_i)-1|\leq \frac{{q}_+^*+{q}_-^*}{2}\leq m_{p_i}f^*.
\end{equation}
Note that $x\frac{\partial H_{n+1}}{\partial x}(x,y)+y\frac{\partial H_{n+1}}{\partial y}(x,y)=(n+1)H_{n+1}(x,y)$. It follows
\begin{equation}\label{hl}
\frac{\partial H_{n+1}}{\partial x}(1,y_i)=\begin{cases}
&-a_1y_i,\quad l=1,\\
&0,\quad l\geq 2.
\end{cases}
\end{equation}
This implies that  $m_{p_i}g^*\geq 1$ if $l\geq 2$. Further
by the properties of Intersection number, we have
$$
m_{p_i} f^*\leq m_{p_i} f^* m_{p_i} g^*\leq I_i.
$$
By \eqref{leq1}, we obtain
$$
i_{P(X)}(p_i)\geq 1-m_{p_i}f^*\geq 1-I_i.
$$
This leads that \eqref{eqi} holds for $l\ge 2$.
On the other hand, if $l=1$, the Jacobian matrix of  \eqref{chartU1} at $p_i$ is
\begin{equation*}
\begin{bmatrix}
 (n+1)a_1 & * \\
 0 &  a_1
\end{bmatrix}.
\end{equation*}
Because $a_1\not=0$,  $p_i$ is an elementary node of \eqref{chartU1}. Thus,   $i_{P(X)}(p_i)=1$. Hence the inequality \eqref{eqi}  holds for $l=1$.

Summarizing the above analysis, we obtain the claim \eqref{eqi} is true. Then we have
\begin{equation}\label{si}
\sum_{inf}i=2\sum_{i=1}^{r} i_{P(X)}(p_i)\geq 2(r-\sum_{i=1}^{r} I_i).
\end{equation}
By  Lemma \ref{lem-f},  Poincar\'{e}-Hopf Theorem and \eqref{si}, we have
\begin{equation}\label{eq-i1}
c(X)+\sum_{i_X(p)\leq 0}i_X(p)=\sum_f i=\frac{1}{2}(2-\sum_{inf}i)\leq 1-r+\sum_{i=1}^r I_i.
\end{equation}
On the other hand, by  Lemma \ref{lem-f}, we have
\begin{equation}\label{eq-i2}
c(X)-\sum_{i_X(p)\leq 0}i_X(p)=\sum_f |i|.
\end{equation}
Adding two inequalities \eqref{eq-i1} and \eqref{eq-i2}, it follows that
$$
c(X)\leq \frac{1-r+\sum_f |i|+\sum_{i=1}^r I_i}{2}.
$$
Note that
\begin{equation*}
\begin{split}
 \sum_f|i|+\sum_{i=1}^r I_i &\leq \sum_{p\in \frac{\partial H}{\partial y}\cap \frac{\partial H}{\partial x}\subset\mathbb{R}^2}I\left(p,\frac{\partial H}{\partial y}\cap \frac{\partial H}{\partial x}\right)+\sum_{i=1}^r I_i\\
&\leq \sum_{p\in\mathbb{CP}^2}I\left(p,\frac{\partial H}{\partial y}\cap \frac{\partial H}{\partial x}\right)=n^2
\end{split}
\end{equation*}
by B\'ezout's Theorem and Lemma \ref{lem2}.  Therefore,  we have $c(X)\leq \frac{n^2+1-r}{2}$. Let
$$
C_{n,m}=C_{n,n}=\left[\frac{n^2+1-r}{2}\right] \ {\rm if}\ n=m.
$$
Since $c(X)$ is a nonnegative integer, $c(X)\leq C_{n,m}$ in the case $n=m$. Hence, we finish the proof of $c(X)\leq C_{n,m}=\left[\frac{n^2+1-r}{2}\right]$ in the case $n=m$.

\textbf{Case 2.} when $n>m$,  $H(x,y)$ has the form
$$
H(x,y)=a_{n+1}y^{n+1}+a_n y^{n}+\dots+a_{m+2}y^{m+2}+H_{m+1}(x,y)+\dots+H_1(x,y).
$$
where $a_{n+1}\neq 0$. Hence, $H_{n+1}(x,y)=a_{n+1}y^{n+1}$ which has a unique linear factor $y$. By Poincar\'e compactification $P(X)$ of $X$, we have the following differential system  in local chart $U_1$
\begin{equation}\label{chartU1-2}
\begin{split}
\frac{dy}{dt}&=-yf^*(1,y,z)+z^{n-m}g^*(1,y,z)\\
&=a_{n+1}y^{n+1}+\sum_{i=1}^n(n+1-i)z^iH_{n+1-i}(1,y),\\
\frac{dz}{dt}&=-zf^*(1,y,z)=z\sum_{i=0}^nz^{i}\frac{\partial H_{n+1-i}}{\partial y}(1,y).
\end{split}
\end{equation}
System \eqref{chartU1-2} has a  unique  critical point $p_0=(0,0)$ which corresponds to one pair of infinite critical points of $X$. We now discuss the estimation on index of system \eqref{chartU1-2} at $p_0$ in two cases:  $g^*(1,0,0)=0$
and $g^*(1,0,0)\neq 0$.

Case (2.i): if $g^*(1,0,0)=0$, that is $\frac{\partial H_{m+1}}{\partial x}(1,0)=0$, then $X\notin \Psi_{n,m}$ by definition of the set $\Psi_{n,m}$. Using the similar arguments on the estimation
of the index at $p_i$ in case 1: $n=m$,  and $l>2$ in \eqref{hl}, we can obtain  the inequality  \eqref{leq1} as follows
$$
|i_{P(X)}(p_0)-1|\leq m_{p_0}f^*.
$$
Note that  $g^*(1,0,0)=0$ implies that $m_{p_0} g^*\geq 1$. Thus we have
$$
i_{P(X)}(p_0)\geq 1-m_{p_0} f^*\geq 1- m_{p_0} f^*m_{p_0}g^*\geq 1-I(p_0,f^*\cap g^*).
$$
By Poincar\'{e}-Hopf Theorem, we have
\begin{equation}\label{eq1-1}
\sum_f i=\frac{1}{2}(2-\sum_{inf} i)\leq I(p_0,f^*\cap g^*).
\end{equation}
By B\'ezout's Theorem and Lemma \ref{lem2}, we have
\begin{equation}\label{eq1-2}
\sum_f|i|+I(p_0,f^*\cap g^*)\leq \sum_{p\in f^*\cap g^*\subset\mathbb{R}^2}I(p,f^*\cap g^*)+I(p_0,f^*\cap g^*)\leq nm.
\end{equation}
Adding two inequalities \eqref{eq1-1} and \eqref{eq1-2},  we know
$$
2c(X)=\sum_f i+\sum_f|i|\leq nm
$$
by Lemma \ref{lem-f}.
Let
$$
C_{n,m}=[\frac{nm+1}{2}]\ {\rm if}\ n>m.
$$
Hence,
 $c(X)\leq \frac{nm}{2}\leq [\frac{nm+1}{2}]=C_{n,m}$ if $n>m$.

Case (2.ii): if $g^*(1,0,0)\neq 0$, that is $\frac{\partial H_{m+1}}{\partial x}(1,0)\not=0$, then $X\in \Psi_{n,m}$.
 Without loss of generality, assume that $g^*(1,0,0)=\frac{\partial H_{m+1}}{\partial x}(1,0)> 0$. Let us estimate the index of $P(X)$ at $p_0$. Taking a circle $S_\epsilon=\{(y,z):y^2+z^2=\epsilon^2,0<\epsilon\ll 1\}$ around $p_0$ and a direction vector $v=(1,0)$ in local chart $U_1$, we consider the points at which the field vector of \eqref{chartU1-2} is parallel to $v$. These points lie on curve $zf^*(1,y,z)=0$, that is,  $z=0$ or $f^*(1,y,z)=0$. Since $g^*(1,0,0)=\frac{\partial H_{m+1}}{\partial x}(1,0)> 0$,  $g^*(1,y,z)>0$ inside $S_\epsilon$ when $\epsilon$ is small enough. Then  ${dy}/{dt}=z^{n-m}g^*(1,y,z)$ on the intersection points of real algebraic curve $f^*(1,y,z)=0$ and $S_\epsilon$, which leads that  ${dy}/{dt}$ has the same sign with $z^{n-m}$. Therefore, near the intersection points of $z=0$ and $S_\epsilon$, the direction vector of $P(X)$ on $S_\epsilon$ can be approximated as
$$
P(X)\approx((n+1)a_{n+1}y^{n+1},(n+1)a_{n+1}y^{n}z).
$$
By analysis the sign of $a_{n+1}$ and the parity of $n$ and $m$, we easily obtain the eight cases on direction vectors $P(X)$ at the intersection points of $zf^*(1,y,z)=0$ and $S_\epsilon$ in counterclockwise sense, see Figure \ref{fig2}. Then we have
\begin{figure}
  \centering
  \includegraphics[width=1\textwidth]{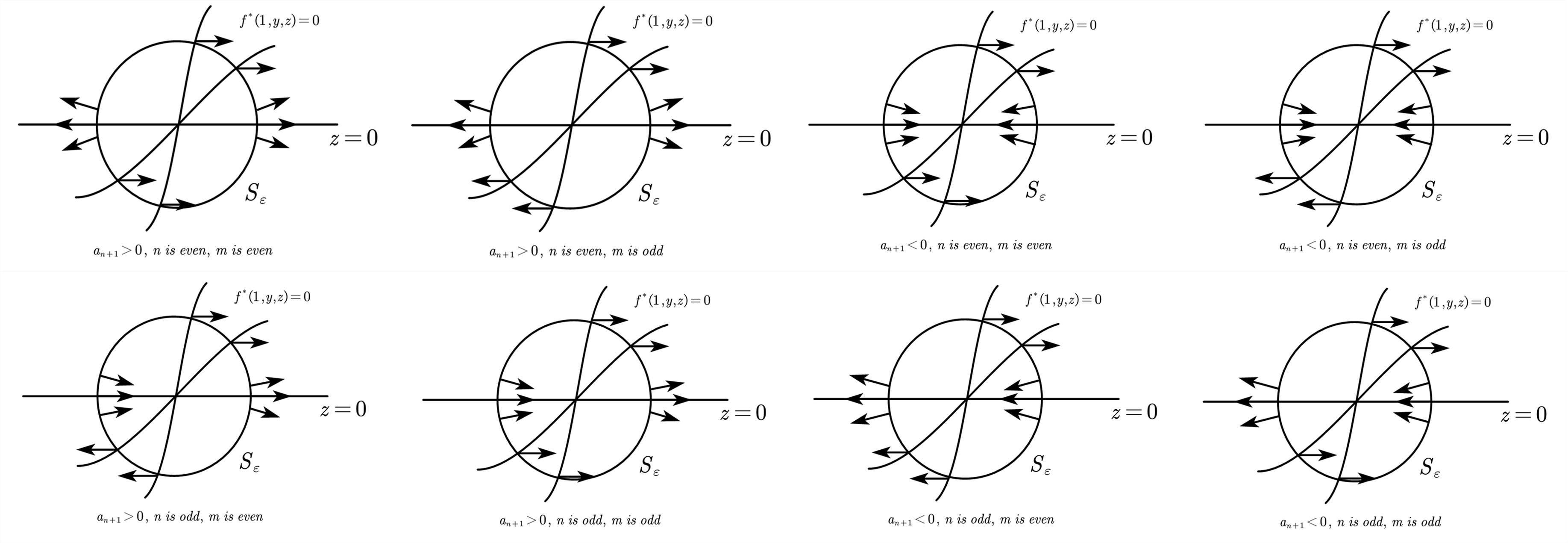}
  \caption{\small The direction vector $P(X)(M)$ on $S_\epsilon$ near the intersection points of $zf^*(1,y,z)=0$ and $S_\epsilon$.}
  \label{fig2}
\end{figure}
\begin{equation}\label{index-inf}
i_{P(X)}(p_0)=\begin{cases}
&0, \quad \text{both $n$ and $m$ are odd and $a_{n+1}>0$,}\\
&+2, \quad \text{both $n$ and $m$ are odd and $a_{n+1}<0$,}\\
&+1, \quad \text{others.}
\end{cases}
\end{equation}
Thus, $\sum_{inf}i=2i_{P(X)}(p_0)\geq 0$. Hence,  we have
$$
c(X)+\sum_{i_X(p)\leq 0}i_X(p)=\sum_f i=\frac{1}{2}(2-\sum_{inf}i)\leq 1
$$
by Poincar\'{e}-Hopf Theorem.
And from B\'ezout's Theorem and Lemma \ref{lem2}, it follows
$$
c(X)-\sum_{i_X(p)\leq 0}i_X(p)=\sum_f |i|\leq nm.
$$
Adding the above two inequalities, we obtain $c(X)\leq \frac{nm+1}{2}$. Notes $c(X)$ is an integer. We have $c(X)\leq [\frac{nm+1}{2}]=C_{n,m}$ if $n>m$.

 Sum up case (2.i) and (2.ii), we finish the proof of $c(X)\leq C_{n,m}=[\frac{nm+1}{2}]$ in the case $n>m$. Hence, the proof of the first part of Theorem \ref{thmmain1} is complete.
\end{proof}

%\begin{corollary}\label{thm-1}
%HK-vector fields of degree $n$ have at most $[\frac{n^2-1}{2}]$ centers on $\mathbb{R}^2$. Furthermore, this bound is the best estimation.
%\end{corollary}

To prove the upper bound $C_{n,m}$ is sharp in Theorem \ref{thmmain1}, we need to construct a Hamiltonian vector field $X\in\mathcal{Y}_{n,m}$ such that $X$ has $C_{n,m}$ centers. Let us first prove  some properties of Hamiltonian vector fields $X$ in $\Psi_{n,n}$, that is, ${\partial H_{n+1}(x,y)}/{\partial y}$ and ${\partial H_{n+1}(x,y)}/{\partial x}$ have no common components.
\begin{lemma}\label{lem-inf-node}
Hamiltonian vector field $X=(-\frac{\partial H}{\partial y},\frac{\partial H}{\partial x})$ is in $\Psi_{n,n}$ if and only if the multiplicity of any a complex polynomial factor of $H_{n+1}(x,y)$ is equal to 1. Furthermore, all infinite critical points of $X$ are elementary nodes if $X\in\Psi_{n,n}$.
\end{lemma}
\begin{proof}
We prove necessity by contradiction. Assume that $H_{n+1}(x,y)$ has polynomial factors with multiplicity $r\ge 2$. Without loss of generality, we assume that there exists $a\in\mathbb{C}$ such that $(y+ax)^r|H_{n+1}(x,y)$. This implies that $-\frac{\partial H_{n+1}(x,y)}{\partial y}$ and $\frac{\partial H_{n+1}(x,y)}{\partial x}$ has a common factor $y+ax$. It is a contradiction with the fact $X\in\Psi_{n,n}$. Hence, the multiplicity of any complex polynomial factor of $H_{n+1}(x,y)$ is equal to 1 if $X=(-\frac{\partial H}{\partial y},\frac{\partial H}{\partial x})\in\Psi_{n,n}$.

 On the other hand, if $X\notin \Psi_{n,n}$, then there exists $(x_0,y_0)\neq (0,0)\in\mathbb{C}^2$ such that $\frac{\partial H_{n+1}}{\partial x}(x_0,y_0)=\frac{\partial H_{n+1}}{\partial y}(x_0,y_0)=0$ by definition of $\Psi_{n,n}$. Note that
$$
(n+1)H_{n+1}(x_0,y_0)=x_0\frac{\partial H_{n+1}}{\partial x}(x_0,y_0)+y_0\frac{\partial H_{n+1}}{\partial y}(x_0,y_0)=0.
$$
Without loss of generality, assume that $x_0=1$. It follows that $H_{n+1}(x,y)$ has a linear factor $y-y_0x$. If the multiplicity of linear factor $y-y_0x$ is equal to 1, then $y=y_0$ is a simple root for the equation $H_{n+1}(1,y)=0$. It is a contradiction to the fact $\frac{\partial H_{n+1}}{\partial y}(x_0,y_0)=0$. Thus,  the proof of sufficiency is finished.

We now consider all infinite critical points of $X$.  Since the arguments are similar in the two local charts $U_1$ and $U_2$, without loss of generality, let all infinite critical points be on a local chart $U_1$ and $P(X)$ has the form \eqref{chartU1} on this chart. $H_{n+1}(1,y_i)=0$ if $(y_i,0)$ is an infinite critical point of $P(X)$ which is equivalent to that $H_{n+1}(x,y)$ has a real linear factor $y-y_ix$.
Therefore, the Jacobian matrix of $P(X)$ at $(y_i,0)$ is given by
\begin{equation*}
\begin{bmatrix}
 (n+1)\frac{\partial H_{n+1}}{\partial y}(1,y_i) & * \\
 0 &  \frac{\partial H_{n+1}}{\partial y}(1,y_i)
\end{bmatrix},
\end{equation*}
which has two different real nonzero eigenvalues $\frac{\partial H_{n+1}}{\partial y}(1,y_i)$ and $(n+1)\frac{\partial H_{n+1}}{\partial y}(1,y_i)$ because the multiplicity of any complex factors of $H_{n+1}(x,y)$ is equal to 1. It follows that all the infinite critical points of $X$ are elementary nodes.
\end{proof}
%\begin{remark}
%\par Notes that following two statement are mutually equivalent:
%\begin{itemize}
%\item[(i)] Multiplicity of any real polynomial factor of $H_{n+1}$ is equal to 1;
%\item[(ii)] Multiplicity of any complex linear factor of $H_{n+1}$ is equal to 1.
%\end{itemize}
%In fact, all the complex polynomial factor of $H_{n+1}$ has linear form $x+ay$ (resp. $y$) if $H_{n+1}(1,a)=0$ (resp. $H(0,1)=0$). If $a\in\mathbb{R}$, $x+ay$ is also a real factor of $H_{n+1}(x,y)$. If $a\in\mathbb{C}\setminus\mathbb{R}$, the complex factor $x+ay$ and $x+\bar{a}y$ appear in pairs and correspond to a real polynomial factor $x^2+(a+\bar{a})x+|a|^2$, where $\bar{a}$ is the complex conjugate of $a$. Thus we only need to prove $X=(-\frac{\partial H}{\partial y},\frac{\partial H}{\partial x})\in\Psi_{n,n}$ is equivalent to that multiplicity of any complex linear factor of $H_{n+1}$ is equal to 1.
%\end{remark}

\begin{lemma}\label{lem-equiv}
Suppose Hamiltonian vector field $X=(-\frac{\partial H}{\partial y},\frac{\partial H}{\partial x})\in\Psi_{n,m}$, and $X$ satisfies one of the following two conditions:
\begin{itemize}
\item[(a)] $m=n$ and $X$ has exactly $2r$ infinite critical points;
\item[(b)]  $n>m$ and $nm$ is even.
\end{itemize}
Then $X$ has at most $nm-C_{n,m}$ saddles,  and the following three statements are equivalent.
\begin{itemize}
\item[(i)] $X$ has $nm$ finite critical points,
\item[(ii)] $X$ has $nm-C_{n,m}$ saddles,
\item[(iii)] $X$ has $C_{n,m}$ centers,
\end{itemize}
where
$ C_{n,m}=\begin{cases}
&[\frac{n^2+1-r}{2}], \quad \text{as}\ n=m,\\
&[\frac{nm+1}{2}], \quad \text{as}\ n>m,
\end{cases}
$
which is defined in Theorem \ref{thmmain1}, and a saddle is an isolated critical point whose neighbourhood is exactly consisting of four hyperbolic sectors.
\end{lemma}

\begin{proof}
First of all, we prove that Hamiltonian vector field $X$ has at most $nm-C_{n,m}$ saddles under the conditions (a) or (b).

If the condition (a) holds, that is $n=m$ and $X$ has exactly $2r$ infinite critical points, then from Lemma \ref{lem-inf-node} we know that $2r$ infinity critical points of $X$ are  elementary nodes. Thus,  $\sum_{inf} i=2r$. By Poincar\'{e}-Hopf Theorem, we have
\begin{equation}\label{a}
\sum_f i=\frac{1}{2}(2-\sum_{inf}i)=\frac{1}{2}(2-2r)=1-r.
\end{equation}
If the condition (b) holds, that is $n>m$ and $nm$ is even, then $X$ has  a unique pair of infinite critical points, without loss of generality, denote them by $p_0=(\pm 1,0,0)$. From the equality \eqref{index-inf}, we have $i_{P(X)}(p_0)=1$. Hence, $\sum_{inf}i=2i_{P(X)}(p_0)=2$. By Poincar\'{e}-Hopf Theorem, we have
\begin{equation}\label{b}
\sum_f i=\frac{1}{2}(2-\sum_{inf}i)=0.
\end{equation}

Let now us denote the number of saddles of  $X$ by $s(X)$. By Lemma \ref{lem-f}  and combining \eqref{a} and \eqref{b}, we have
\begin{equation}\label{eq2-3}
s(X)-c(X)\leq -\sum_f i=\begin{cases}
&r-1, \quad  \text{as}\ n=m,\\
&0, \quad \text{as}\ n>m \text{ and $nm$ is even},
\end{cases}
\end{equation}
On the other hand, by B\'ezout's Theorem, we have
\begin{equation}\label{eq2-4}
s(X)+c(X)\leq nm.
\end{equation}
Adding the inequality \eqref{eq2-3} and \eqref{eq2-4}, we obtain that
$$
s(X)\leq \begin{cases}
&\frac{n^2+r-1}{2}=n^2-\frac{n^2+1-r}{2}, \quad  \text{as}\ n=m,\\
&\frac{nm}{2}=nm-\frac{nm}{2}, \quad \text{as}\ n>m \text{ and $nm$ is even}.
\end{cases}
$$
Note that $r$ pairs of infinity critical points of $X$ correspond to all different real linear factors of $H_{n+1}(x,y)$. Because the  multiplicity of any real factors of $H_{n+1}(x,y)$ is equal to 1 by Lemma \ref{lem-inf-node},  $n+1-r$ is even, that is $n\not\equiv r$ (mod 2). Thus $n^2+1-r$ is even which implies that $\frac{n^2+1-r}{2}=[\frac{n^2+1-r}{2}]=C_{n,n}$. This leads that $s(X)\leq n^2-C_{n,n}$ when $n=m$. It is clear that $\frac{nm}{2}=[\frac{nm+1}{2}]$ when $nm$ is even, which leads that $s(X)\leq nm-C_{n,m}$ when $n>m$ and $nm$ is even.

Next we prove statements $(i),(ii)$ and $(iii)$ are equivalent.

 (i) $\Rightarrow$ (ii):
 Since Hamiltonian  vector field $X$ has $nm$ finite critical points in $\mathbb{R}^2$, each a critical point is one of elementary center and elementary saddle by the property of the intersection number. Then
 $$
 c(X)+s(X)=nm \ \Rightarrow\ \  c(X)=nm-s(X).
 $$
If $n=m$, we have
$$
(n^2-s(X))-s(X)=\sum_f i=1-r.
$$
Hence, $s(X)=\frac{n^2+r-1}{2}=n^2-\frac{n^2+1-r}{2}$. Note $\frac{n^2+1-r}{2}=[\frac{n^2+1-r}{2}]=C_{n,n}$ when $X$ belongs to $\Psi_{n,n}$. This leads that $s(X)=n^2-C_{n,n}$.

If $n>m$ and $nm$ is even, we have
$$
(nm-s(X))-s(X)=\sum_f i=0.
$$
Hence, $s(X)=\frac{nm}{2}=nm-\frac{nm}{2}$. Note $\frac{nm}{2}=[\frac{nm+1}{2}]=C_{n,m}$ when $nm$ is even. This leads that $s(X)=nm-C_{n,m}$.

(ii) $\Rightarrow$ (iii): By Poincar\'{e}-Hopf Theorem and Lemma \ref{lem-f}, we have
$$
c(X)-(nm-C_{n,m})\geq\sum_f i=\begin{cases}
&1-r, \quad  \text{as}\ n=m,\\
&0, \quad \text{as}\ n>m \text{ and $nm$ is even},
\end{cases}
$$
since $X$ has $nm-C_{n,m}$ saddles.
It follows that
$$
c(X)\geq nm+\sum_f i-C_{n,m}=\begin{cases}
&\frac{n^2+1-r}{2}, \quad \text{as}\ n=m,\\
&\frac{nm}{2}, \quad \text{as}\ n>m \text{ and $nm$ is even},
\end{cases}
$$
i.e. $c(X)\geq C_{n,m}$. On the other hand, $c(X)\le C_{n,m}$ by the first part of Theorem \ref{thmmain1}. Therefore, $c(X)=C_{n,m}$.

(iii) $\Rightarrow$ (i): Assume that $p_1,p_2,\dots,p_l$ are finite critical points of $X$ which are not centers. Then
$$
i_X(p_i)\leq 0,\quad \forall i=1,2,\dots, l
$$
by Lemma \ref{lem-f}. Note that $\sum_{inf} i=2r$. And when $X=(-\frac{\partial H}{\partial y},\frac{\partial H}{\partial x})\in\Psi_{n,m}$, we recall that
$$C_{n,m}=\begin{cases}
&[\frac{n^2+1-r}{2}]=\frac{n^2+1-r}{2}, \quad  \text{as}\ n=m,\\
&[\frac{nm+1}{2}]=\frac{nm}{2}, \quad \text{as}\ n>m \text{ and $nm$ is even}.
\end{cases}  $$
By Poincar\'{e}-Hopf Theorem, we have
\begin{equation}\label{se}
\sum_{i=1}^l i_X(p_i)=\sum_f i-C_{n,m}=\begin{cases}
&1-r-\frac{n^2+1-r}{2}=-\frac{n^2+r-1}{2}, \quad  \text{as}\ n=m,\\
&-\frac{nm}{2}, \quad \text{as}\ n>m \text{ and $nm$ is even}.
\end{cases}
\end{equation}

Let $I_i=I(p_i,(-\frac{\partial H}{\partial y})\cap\frac{\partial H}{\partial x})$. By B\'ezout's Theorem,
$$
C_{n,m}+\sum_{i=1}^l I_i\leq \sum_{p\in(-\frac{\partial H}{\partial y})\cap\frac{\partial H}{\partial x}\subseteq\mathbb{R}^2}I(p,(-\frac{\partial H}{\partial y})\cap\frac{\partial H}{\partial x})\leq\sum_{p\in\mathbb{CP}^2}I(p,(-\frac{\partial H}{\partial y})\cap\frac{\partial H}{\partial x})=nm.
$$
By using Lemma \ref{lem2}, we have
\begin{equation*}
\begin{split}
-\sum_{i=1}^l i_X(p_i)&=\sum_{i=1}^l |i_X(p_i)|\leq \sum_{i=1}^l\sqrt{I_i}\leq \sum_{i=1}^l I_i\\
& \leq nm-C_{n,m}=\begin{cases}
&\frac{n^2+r-1}{2}, \quad \text{as}\ n=m,\\
&\frac{nm}{2}, \quad \text{as}\ n>m \text{ and $nm$ is even},
\end{cases}
\end{split}
\end{equation*}
which implies that
$$
\sum_{i=1}^l |i_X(p_i)|=\sum_{i=1}^l\sqrt{I_i}=\sum_{i=1}^l I_i=nm-C_{n,m}
$$
by \eqref{se}. Therefore,
 $$
 I_i=1, \ i_X(p_i)=-1, \ l=nm-C_{n,m}, \ i=1, \cdots, l,
 $$
 because $I_i$ is positive integer and $i_X(p_i)$ is nonpositive integer.
 Hence, the number of finite critical points of $X$ is
 $$C_{n,m}+l=C_{n,m}+(nm-C_{n,m})=nm.$$

 This lemma is proved.
\end{proof}

%\begin{remark}
{\bf Remark 3.1} Lemma \ref{lem-equiv} tell us that $X\in G_{n,n}$ if $X$ has $C_{n,n}$ centers. This improves Proposition 3.5 in \cite{CGM}. However,
if both $n$ and $m$ are odd with $n>m$,  it can be proved that the statements  $(iii) \Rightarrow (ii)$ and $(iii) \Rightarrow (i)$ are true by following the proof of Lemma \ref{lem-equiv}. But $(i) \nRightarrow (ii)$ and $(i) \nRightarrow (iii)$, please see the following system
 \begin{equation}\label{example}
\begin{split}
\dot{x}&=y(y-1)\dots (y-(n-1)),\\
\dot{y}&=x(x-1)\dots (x-(m-1)),
\end{split}
\end{equation}
where both $n$ and $m$ are odd with $n>m$. It is easily checked that system \eqref{example} has $nm$ finite critical points, in which there are exactly $\frac{nm-1}{2}$ centers and $\frac{nm+1}{2}$ saddles. Thus, $(i) \nRightarrow (ii)$ and $(i) \nRightarrow (iii)$. Moreover, $$\frac{nm+1}{2}>nm-[\frac{nm+1}{2}]=nm-C_{n,m},$$ which implies that the number of saddles of system \eqref{example} is greater than $nm-C_{n,m}$. This leads that the conclusion in Lemma \ref{lem-equiv}, $X$ has  at most $nm-C_{n,m}$ saddles, does not hold if $n>m$ and both $n$ and $m$ are odd.
%\end{remark}

\medskip

We are now in the position to prove  the second conclusion of Theorem \ref{thmmain1}: the upper bound $C_{m,n}$ is sharp.

\begin{proof}[Proof of the second part of Theorem \ref{thmmain1}]
We shall construct a polynomial Hamiltonian vector field $X\in\mathcal{Y}_{n,n}$ which has $C_{n,m}$ centers. When $n>m$, an example has been given in \cite{CGM}. So we only consider the case $n=m$, and distinguish two cases depending on the relation between $n$ and $r$.
\par \textbf{Case I.} $n \not\equiv r$ (mod 2).
\par In this case, both $\frac{n+1-r}{2}$ and $\frac{n^2+1-r}{2}$ are  integers. Let
\begin{equation}\label{Hamifun}
H(x,y)=\prod_{i=1}^{\frac{n+r+1}{2}} \bar{H}_i
\end{equation}
where
\begin{equation*}
\bar{H}_i(x,y)=\left\{
\begin{aligned}
&4^ix^2+\frac{y^2}{4^i}-1,\quad \text{as}\ 1\leq i\leq \frac{n+1-r}{2},\\
&2^ix+\frac{y}{2^i}-1,\quad \text{as}\ \frac{n+1-r}{2}+1\leq i\leq \frac{n+1+r}{2}.
\end{aligned}
\right.
\end{equation*}
Then $$H_{n+1}(x,y)=\prod_{i=1}^{\frac{n+1-r}{2}}\left(4^ix^2+\frac{y^2}{4^i}\right)\prod_{i=\frac{n+1-r}{2}+1}^{\frac{n+r+1}{2}}\left(2^ix+\frac{y}{2^i}\right).$$
Thus, the Hamilton vector field $X$ with Hamiltonian function $H(x,y)$ in \eqref{Hamifun} are in $\Psi_{n,n}$ by Lemma \ref{lem-inf-node}.
It is easily checked that the following conclusions are true.
\begin{itemize}
\item[(A)] $\sharp\{(x,y)\in\mathbb{R}^2: \bar{H}_i=\bar{H}_j=0, i\neq j\}=deg (\bar{H}_i)\cdot deg (\bar{H}_j)$;
\item[(B)] $\{(x,y)\in\mathbb{R}^2: \bar{H}_i=\bar{H}_j=\bar{H}_k=0, i\neq j, j\neq k, i\neq k\}=\emptyset$;
\item[(C)] $H_{n+1}(x,y)$ has exactly $r$ different real linear factors $2^ix+\frac{y}{2^i}$, $i=\frac{n+1-r}{2}+1, \cdots, \frac{n+1+r}{2}$.
\end{itemize}

To show that $X$ with Hamiltonian function \eqref{Hamifun} has $\frac{n^2+1-r}{2}$ centers, it is enough to prove that $X$ have $\frac{n^2+r-1}{2}$ saddles by Lemma \ref{lem-equiv}. Let us to study the critical points of algebraic curve $H(x,y)=0$  in $\mathbb{R}^2$. It is clear that $H(x,y)$ in \eqref{Hamifun} has $r$ linear factors $2^ix+\frac{y}{2^i}-1$ with $\frac{n+1-r}{2}+1\leq i\leq \frac{n+1+r}{2}$,  and $\frac{n+1-r}{2}$ quadratic polynomial factors $4^ix^2+\frac{y^2}{4^i}-1$ with $1\leq i\leq \frac{n+1-r}{2}$. Therefore,  the common points of either two lines or two ellipses or one ellipse and one line, that is the points  $(x_0,y_0)\in \{\bar{H}_k=0\}\cap \{\bar{H}_l=0\}$ with $k\not=l$,  are critical points of $X$. By straightforward calculation, we obtain that
\begin{itemize}
\item[(i)] the common points of two lines have  $\frac{r(r-1)}{2}$;
\item[(ii)] the common points of two ellipses have  $\frac{1}{2}(n+1-r)(n-1-r)$;
\item[(iii)] the common points of one ellipse and one line have $r(n+1-r)$.
\end{itemize}
Hence, the total common points are
$$
\frac{r(r-1)}{2}+\frac{1}{2}(n+1-r)(n-1-r)+r(n+1-r)=\frac{n^2+r-1}{2}.
$$
 It is easily checked that those common points are saddles of $X$ with Hamiltonian function \eqref{Hamifun} since the determinant of Jacobian matrixes of $X$ at any a common point $(x_0,y_0)\in \{\bar{H}_k=0\}\cap \{\bar{H}_l=0\}$ is
 $$
-\left(\prod_{i=1,i\neq k,l}^{\frac{n+r+1}{2}}\bar{H}_i(x_0,y_0)\right)^2\cdot\left(\begin{vmatrix}
\frac{\partial \bar{H}_k}{\partial x} & \frac{\partial \bar{H}_k}{\partial y} \\
\frac{\partial \bar{H}_l}{\partial x} & \frac{\partial \bar{H}_l}{\partial y}
\end{vmatrix}{(x_0,y_0)}\right)^2<0
$$
by (A) and (B). Therefore, $s(X)\ge \frac{n^2+r-1}{2}$. By Lemma \ref{lem-equiv} we know that $s(X)\le \frac{n^2+r-1}{2}$.  Hence,  $s(X)=\frac{n^2+r-1}{2}$. This implies that $X$ has $\frac{n^2+1-r}{2}$ centers as $n \not\equiv r$ (mod 2).

 \textbf{Case II.} $n\equiv r$ (mod 2). In this case, $n+1-r$ is odd, so Hamiltonian vector field $X\in\mathcal{Y}_{n,n}\setminus\Psi_{n,n}$, which is non-generic. Let us construct a Hamiltonian function %$r\geq 1$ since $H_{n+1}(x,y)$ has at least one linear factor if $n$ is even.
\begin{equation}\label{Hamifun2}
H(x,y)=\prod_{i=1}^{\frac{n+r}{2}} \hat{H}_i
\end{equation}
where
\begin{equation*}
\hat{H}_i=\left\{
\begin{aligned}
&x-\frac{1}{2}y^2+2,\quad i=1,\\
&4^ix^2+\frac{y^2}{4^i}-1,\quad 2\leq i\leq \frac{n-r}{2}+1,\\
&2^ix+\frac{y}{2^i}-1,\quad \frac{n-r}{2}+2\leq i\leq \frac{n+r}{2}.
\end{aligned}
\right.
\end{equation*}
Then
$$
H_{n+1}(x,y)=-\frac{1}{2}y^2\prod_{i=2}^{\frac{n-r}{2}+1}\left(4^ix^2+\frac{y^2}{4^i}\right)\prod_{i=\frac{n-r}{2}+2}^{\frac{n+r}{2}}\left(2^ix+\frac{y}{2^i}\right).
$$
It can be checked that conclusions (A) and (B) in Case I still hold,  and $H_{n+1}(x,y)$ has exactly $r$ different linear factors $y$ and $2^ix+\frac{y}{2^i}$, where $\frac{n-r}{2}+2\leq i\leq \frac{n+r}{2}$. Using the similar arguments in Case I, we can prove that the Hamiltonian vector field $X$ with Hamiltonian function \eqref{Hamifun2} has  $\frac{n^2+r-2}{2}$ finite critical points which are all elementary saddles.

To obtain the number $c(X)$ of centers of $X$ with Hamiltonian function \eqref{Hamifun2}, we  claim that $\sum_{inf} i=2r$ for $X$. If it is true, by Lemma \ref{lem-f} and \eqref{phT}, we have
$$
c(X)-\frac{n^2+r-2}{2}\geq \sum_f i=\frac{1}{2}(2-\sum_{inf} i)=1-r.
$$
It follows $c(X)\geq \frac{n^2-r}{2}=[\frac{n^2+1-r}{2}]$ when $n\equiv r$ (mod 2). From the proof of first part of Theorem \ref{thmmain1}, we have $c(X)=[\frac{n^2+1-r}{2}]$, that is $X$ with Hamiltonian function \eqref{Hamifun2} has $[\frac{n^2+1-r}{2}]$ centers. Then the proof is completed.

 Now we prove the claim $\sum_{inf} i= 2r$. It is to calculate the index of every infinite critical points of $X$. Note that all infinite critical points of $X$ come from the lines $2^ix+\frac{y}{2^i}=0$ with $\frac{n-r}{2}+2\leq i\leq \frac{n+r}{2}$ and a parabola $y^2$ by the expression of $H_{n+1}(x,y)$.

 Using the similar arguments in the proof of Lemma \ref{lem-inf-node}, we can see that $r-1$ pairs of infinite critical points  corresponding to the lines $2^ix+\frac{y}{2^i}=0$ with $\frac{n-r}{2}+2\leq i\leq \frac{n+r}{2}$ are all elementary nodes and each of their indexes is $+1$.
The only issue left is to calculate the index at one pair of infinite critical points  corresponding to $y^2=0$. The method used is the  same to that in the proof of first part of Theorem \ref{thmmain1}. Recall the Poincar\'e compactification $P(X)$ of $X$ has the form \eqref{chartU1} on the local chart $U_1$. Taking the circle $S_\epsilon=\{(y,z):y^2+z^2=\epsilon^2,0<\epsilon\ll 1\}$ and a direction vector $v=(1,0)$ on the chart $U_1$, we consider the points on $S_\epsilon$ at which the field vector $P(X)$ is parallel to $v$, i.e. $(-\epsilon,0),(\epsilon,0)$ and the intersection points of $S_\epsilon$ and curve $f^*(1,y,z)=0$. Denote $\prod_{i=2}^{\frac{n+r}{2}} \hat{H}_i(x,y)$ by $\bar{H}(x,y)$. The vector field of $P(X)$ near points $(-\epsilon,0)$ and $(\epsilon,0)$ can be approximate as
$$
P(X)\approx(-\frac{(n+1)}{2}\bar{H}^*(1,0,0)y^2,-\bar{H}^*(1,0,0)yz),
$$
where $\bar{H}^*(x,y,z)$ is the homogenization of polynomial $\bar{H}(x,y)$.
It is easily calculated that
$$
\bar{H}^*(1,0,0)=\prod_{i=2}^{\frac{n-r}{2}+1} 4^i\cdot\prod_{i=\frac{n-r}{2}+2}^{\frac{n+r}{2}} 2^i>0.
$$
Further we have
$$
(\frac{\partial \bar{H}}{\partial x})^*(1,0,0)>0,\quad (\frac{\partial \bar{H}}{\partial y})^*(1,0,0)>0.
$$
Then we calculate the contribution of points $(-\epsilon,0)$ and $(\epsilon,0)$ to index at infinite critical point $(0,0)$. Consider the curve $f^*(1,y,z)=0$, it can be written as
$$
f^*(1,y,z)=y\bar{H}^*(1,y,z)-(z-\frac{1}{2}y^2+2z^2)(\frac{\partial \bar{H}}{\partial y})^*(1,y,z)=0.
$$
Then we have
$$
\nabla f^*(1,y,z)=(\frac{\partial f^*}{\partial y}(1,y,z),\frac{\partial f^*}{\partial z}(1,y,z))=(\bar{H}^*(1,0,0),-(\frac{\partial \bar{H}}{\partial y})^*(1,0,0)),
$$
which means $m_{(0,0)}f^*=1$. Hence, by Implicit Function Theorem, $f^*(1,y,z)=0$ can be regarded as a one-dimension manifold locally near the point $(0,0)$. Notes the tangent line of $f^*(1,y,z)=0$ at $(0,0)$ is $\bar{H}^*(1,0,0)y-(\frac{\partial \bar{H}}{\partial y})^*(1,0,0)z=0$. So when $\epsilon$ is small enough, curve $f^*(1,y,z)=0$ has two intersection points $(y_1,z_1),(y_2,z_2)$ with $S_\epsilon$ such that $y_1>0,z_1>0$ and $y_2<0,z_2<0$.
\begin{figure}
  \centering
  \includegraphics[width=0.5\textwidth]{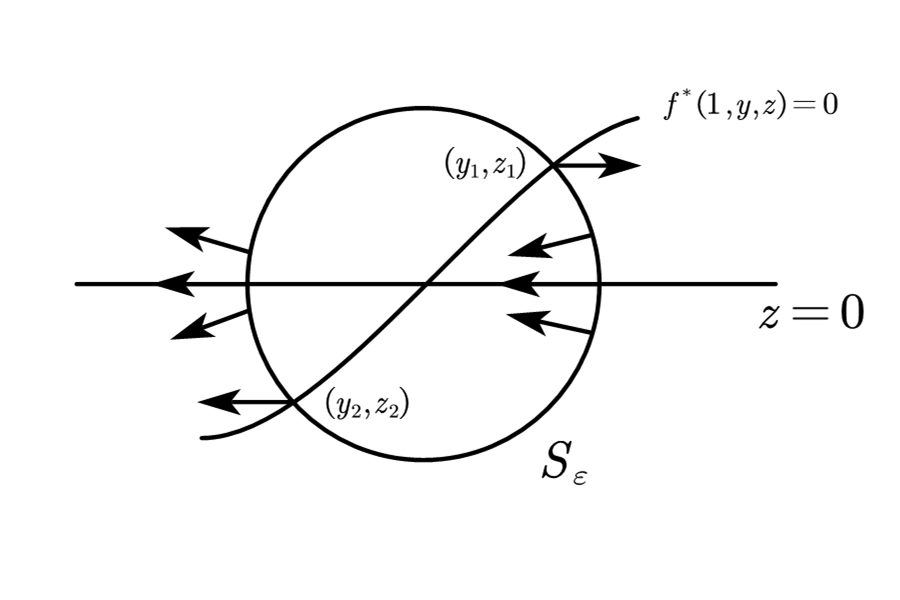}
  \caption{The director vector of $P(X)$ near the intersection points of  $y=0$ and $S_\epsilon$.}
  \label{fig3}
\end{figure}
Near the intersection points of curve $f^*(1,y,z)=0$ and $S_{\epsilon}$, let us calculate the $y$ direction of vector field $P(X)$
\begin{equation*}
\begin{split}
\dot{y}&=g^*(1,y,z)\\
&=z\bar{H}^*(1,y,z)+(z-\frac{1}{2}y^2+2z^2)(\frac{\partial \bar{H}}{\partial x})^*(1,y,z)\\
&=z\bar{H}^*(1,y,z)+\frac{y\bar{H}^*(1,y,z)}{(\frac{\partial \bar{H}}{\partial y})^*(1,y,z)}(\frac{\partial \bar{H}}{\partial x})^*(1,y,z)\\
&=\frac{\bar{H}^*(1,y,z)}{(\frac{\partial \bar{H}}{\partial y})^*(1,y,z)}\left(z(\frac{\partial \bar{H}}{\partial y})^*(1,y,z)+y(\frac{\partial \bar{H}}{\partial x})^*(1,y,z)\right).
\end{split}
\end{equation*}
Since $\bar{H}^*(1,0,0)>0$, $(\frac{\partial \bar{H}}{\partial x})^*(1,0,0)>0$ and $(\frac{\partial \bar{H}}{\partial y})^*(1,0,0)>0$,  we have
$$
\dot{y}|_{(y,z)=(y_1,z_1)}>0,\quad \dot{y}|_{(y,z)=(y_2,z_2)}<0.
$$
So the direction vector  of $P(X)$ near these intersection points can be sketched in  Figure \ref{fig3}. Hence,  the index at the infinite critical point $(0,0)$ in the chart $U_1$ is $+1$ too.

Summarizing the above analysis, we have proved all indices at  the $r$ pairs of infinite critical points are  +1. Thus $\sum_{inf} i=2r$ and the claim is proved. So far, we finish the proof of Theorem \ref{thmmain1}.
\end{proof}

{\bf Remark 3.2}:  Theorem 3.1 in \cite{CGM} is our result (Theorem 3.1) when $r=0$ in the case $m=n$. From the construction of polynomial Hamiltonian vector fields $X$ having $C_{n,m}$ centers in proof of our Theorem 3.1, we can see that the existence of invariant straight lines
reduces the number of centers. If the $r$ pairs of infinite critical points of $X$ come from $r$ invariant straight lines, then the maximum number of centers reduces to $[\frac{n^2+1-r}{2}]$. This implies that the number of real linear factors of polynomial Hamiltonian functions $H(x,y)$ affects the number of ovals of level sets $H(x,y)=h$ in $\mathbb{R}^2$. And topological classifications of infinite critical points of $X$ play a key role in understanding the geometry of real algebraic curve $H(x,y)$ in $\mathbb{R}^2$. If  the number of center of $X$ arrives the least upper bound $C_{n,m}$, then all infinite critical points of $X$ are elementary nodes or singularities with index $+1$, respectively. If $X$ has a unique finite critical point, and all infinite critical points of $X$ have exactly two hyperbolic sectors, then every level sets $H(x,y)=h$ are ovals in $\mathbb{R}^2$  (cf. \cite{HLX}).

\smallskip

It is well known that there are three types center: elementary center, nilpotent center and degenerate center for vector fields (cf. \cite{HLX,LX}).  In the following we characterize the type of center if Hamiltonian vector fields have  $C_{n,m}$ centers.
\begin{proposition}\label{pro-center}
Suppose that Hamiltonian vector field $X=(-\frac{\partial H}{\partial y},\frac{\partial H}{\partial x})\in\mathcal{Y}_{n,m}$ with $2r$ infinite critical points has $C_{n,m}$ centers, where $r$ is a nonnegative integer. Then all centers of $X$ are elementary.
\end{proposition}
\begin{proof}
 We first prove that for any a polynomial Hamiltonian vector field $X=(-\frac{\partial H}{\partial y},\frac{\partial H}{\partial x})\in\mathcal{Y}_{n,m}$ with $2r$ infinite critical points,  there exists a family of polynomial Hamiltonian  vector fields $X_\epsilon$ with at least $2r$ infinite critical points for any $0<\epsilon\ll 1$ such that $X_\epsilon\in\Psi_{n,m}$ and $\lim_{\epsilon \to 0}H(x,y,\epsilon) = H(x,y)$.
 It is obvious if $X\in\Psi_{n,m}$. We only need to study $X\in\mathcal{Y}_{n,m}\setminus\Psi_{n,m}$.

  Consider a perturbation of Hamiltonian function $H(x,y)$ if $X\in\mathcal{Y}_{n,m}\setminus\Psi_{n,m}$, we distinguish two cases: $n=m$ and $n>m$.

 When $n=m$, since $X$ has $2r$ infinite critical points,  $H_{n+1}(x,y)$ has exactly $r$ different linear factors. Thus $H_{n+1}(x,y)$ has the following expression
$$
H_{n+1}(x,y)=\prod_{i=1}^rL_i(x,y)^{k_i}\prod_{i=1}^l M_i(x,y)^{l_i},
$$
where $k_i,l_i\geq 1$ and $\sum_{i=1}^r k_i+\sum_{i=1}^l l_i=n+1$,  $L_1, \cdots, L_r$ are the $r$ different real homogeneous linear polynomials and $M_1,\cdots, M_l$ are the $l$ different irreducible  real homogeneous quadratic polynomials. Let
$$
\tilde{L}_i(x,y)=\prod_{j=1}^{k_i} (L_i(x,y)+\frac{\epsilon}{j} x),\quad \tilde{M}_i(x,y)=\prod_{j=1}^{l_i} (M_i(x,y)+\frac{\epsilon}{j} x^2), \ 0<\epsilon\ll 1
$$
and
$$
H_{n+1}(x,y,\epsilon)=\prod_{i=1}^r\tilde{L}_i(x,y)\prod_{i=1}^l\tilde{M}_i(x,y).
$$
Then the multiplicity of every factors of $H_{n+1}(x,y,\epsilon)$ is one and $H_{n+1}(x,y,\epsilon)$ has $\sum_{i=1}^r k_i$ different real linear factors and $\sum_{i=1}^r k_i\geq r$. It can be checked that
$$\lim_{\epsilon\to 0}H_{n+1}(x,y,\epsilon) = H_{n+1}(x,y).$$
 Let
$$
H(x,y,\epsilon)=H(x,y)-H_{n+1}(x,y)+H_{n+1}(x,y,\epsilon).
$$
Then $\lim_{\epsilon\to 0}H(x,y,\epsilon) = H(x,y)$, the Hamiltonian  vector fields $X_\epsilon$ with Hamiltonian function $H(x,y,\epsilon)$ have at least $2r$ infinite critical points for any $0<\epsilon\ll 1$,  and $X_\epsilon\in\Psi_{n,n}$ by Lemma \ref{lem-inf-node}.

When $n>m$, $X\in\mathcal{Y}_{n,m}\setminus\Psi_{n,m}$ leads to $\frac{\partial H_{m+1}}{\partial x}(1,0)=0$. Let
$$
H(x,y,\epsilon)=H(x,y)+\epsilon x^{m+1},\ 0<\epsilon\ll 1.
$$
Then $\frac{\partial H_{m+1}(x,y,\epsilon)}{\partial x}(1,0,\epsilon)=\epsilon>0$, that is, $\frac{\partial H_{m+1}(x,y,\epsilon)}{\partial x}$ has no factors $y$. This implies that $\frac{\partial H_{m+1}(x,y,\epsilon)}{\partial x}$ and $\frac{\partial H_{n+1}(x,y,\epsilon)}{\partial y}$ do not have common factors. Denote Hamiltonian vector fields $(-\frac{\partial H(x,y,\epsilon)}{\partial y},\frac{\partial H(x,y,\epsilon)}{\partial x})$ by $X_\epsilon$. Hence, $X_\epsilon\in\Psi_{n,m}$, $\lim_{\epsilon\to 0}H(x,y,\epsilon) = H(x,y)$ and $X_\epsilon$  has only two infinite critical points, that is $r=1$.

We then prove that $X_\epsilon$ has exactly $nm$ finite critical points in $\mathbb{R}^2$.

Let $p$ be a finite critical point of $X$. By the additivity of the indices of vector fields (see \cite{Ar3} for detail), the index at $p$ of $X$, which is equal to $+1$, is the sum of indices at those critical points of vector field $X_\epsilon$ which tend to $p_0$ when $\epsilon\rightarrow 0$. Notes all the indices at critical points of $X_\epsilon$ is less than or equal to $+1$ by Lemma \ref{lem-f}. This argument implies vector field $X_\epsilon$ have at least $C_{n,m}$ centers. Denoted the least upper bound of centers of $X_\epsilon$
by $C_{n,m}^{\epsilon}$. Then $C_{n,m}^{\epsilon}\geq C_{n,m}$.  On the other hand, by Theorem \ref{thmmain1}, $X_\epsilon$ has at most $C_{n,m}^{\epsilon}$ centers, where
$$
C_{n,m}^{\epsilon}=\begin{cases}
&[\frac{n^2+1-\sum_{i=1}^r r_i}{2}]\le [\frac{n^2+1-r}{2}]=C_{n,n}, \quad \text{as}\  n=m,\\
&[\frac{nm+1}{2}]=C_{n,m}, \quad \text{as}\ n>m.
\end{cases}
$$
Hence, we have $C_{n,m}^{\epsilon}=C_{n,m}$ and $X_\epsilon$ has exactly $C_{n,m}^{\epsilon}$ centers. By Lemma \ref{lem-equiv} and Remark 3.1, $X_\epsilon$ has exactly $nm$ finite critical points and all the finite critical points of $X_\epsilon$ are elementary.

Finally, we prove the conclusion, all centers of $X$ are elementary, by contradiction. Assume that $X$ has a non-elementary center $p_0$ whose index is $+1$. Then the intersection number, denoted by $i_{p_0}$, of $H_x(x, y)$ and $H_y(x, y)$ at $p_0$ is larger than $1$. Thus when $\epsilon$ is sufficiently small,  the sum of  the intersection numbers of  $H_x(x, y, \epsilon)$  and $H_y(x, y, \epsilon)$ at the intersection points on a neighborhood $U_{p_0}$ of
$p_0$ in $\mathbb C^2$,
is also $i_{p_0}$. Since all the critical points of $X_\epsilon$
are real, the intersection points of   $H_x(x, y, \epsilon)$ and $H_y(x, y, \epsilon)$ on $U_{p_0}$ are real. At last, when $\epsilon$ is sufficiently small, the sum of the indices of
critical points of $X_\epsilon$ on $U_{p_0}$ is the index of $X$ at $p_0$, which is $+1$. Note that the indices at the finite critical points of $X_\epsilon$ are either $+1$ or $-1$. Hence, we must have that there are two centers of $X_\epsilon$ on $U_{p_0}$, which shows that $X_\epsilon$ has at least $C_{n,m}+1$ centers. It is a contradiction with the fact that $X_\epsilon$ has exactly $C_{n,m}$ centers. This finishes the proof.
\end{proof}

\section{Configurations of centers of  Hamiltonian Kolmogorov  systems}\label{configuration}

In this section we study Hamiltonian polynomial vector fields having two intersecting invariant straight lines in $\mathbb{R}^2$.  Since with an affine transformation these two invariant straight lines
can become the axes of coordinates, we consider vector fields $X_{hk}$ and investigate the possible configurations of centers  when $X_{hk}$   has $C_{n,m}$ centers. The corresponding differential system of $X_{hk}$ is
\begin{equation}\label{xhk}
\begin{split}
\frac{dx}{dt}=&-x\left(F(x,y)+y\frac{\partial F}{\partial y}\right),\\
\frac{dy}{dt}=& y\left(F(x,y)+x\frac{\partial F}{\partial x}\right),
\end{split}
\end{equation}
where $F(x,y)$ is a polynomial of degree $n-1$ and $H(x,y)=xyF(x,y)$. It is clear that system \eqref{xhk} has two intersecting invariant straight lines $x=0$ and $y=0$. Hence, the center of system \eqref{xhk} can only be in the interior of four quadrants of $\mathbb{R}^2$.
We say system \eqref{xhk} has {\it a configuration $(i_1;i_2;i_3;i_4)$ of centers if there exist exactly $i_j$ centers in the interior of the $j$th quadrant of $\mathbb{R}^2$ for  $j=1,2,3,4$, where $i_j$ is a nonnegative integer}. Note that linear transformations do not change the total number of centers of system  \eqref{xhk}. Using linear transformations of  the variables $x$ and $y$ if necessary, we can  assume that $i_1=\min\{i_1,i_2,i_3,i_4\}$ and $i_2\leq i_4$ since  the configurations $(i_2;i_1;i_4;i_3),(i_4;i_3;i_2;i_1)$ and $(i_1;i_4;i_3;i_2)$ of centers  can be obtained by the following transformations
\begin{equation}\label{trans}
(x,y)\mapsto(-x,y),\ (x,y)\mapsto(x,-y),\ (x,y)\mapsto (y,x),
\end{equation}
respectively if system \eqref{xhk} has configuration $(i_1;i_2;i_3;i_4)$ of centers. Therefore,  we say {\it two configurations of centers are equivalent} if there exist some transformations in \eqref{trans} such that one configuration of centers can be transformed to the other by these transformations or their composition. Hereafter we consider  all different configurations of centers  of system \eqref{xhk}  in this equivalent sense if the number of centers of system \eqref{xhk} reaches the least upper bound $C_{n,m}$.

From Theorem \ref{thmmain1} and Proposition \ref{pro-center}, we obtain the least upper bound of centers of system \eqref{xhk} and property of centers as follows.
\begin{proposition}\label{thm-HK}
If system \eqref{xhk} is a HK polynomial system of degree $n$ in $\mathbb{R}^2$, then system \eqref{xhk} has at most $[\frac{n^2-1}{2}]$ centers. Furthermore,  this bound $[\frac{n^2-1}{2}]$ is sharp,  and all the centers are elementary if system \eqref{xhk} has $[\frac{n^2-1}{2}]$ centers.
\end{proposition}
By Proposition \ref{thm-HK}, we can see that the number of centers depends on the degree $n$ of this system, and the dynamics of system \eqref{xhk} are trivial when $n=1,2$. If $n=1$, then system \eqref{xhk} is a linear system which has no center. If $n=2$, then system \eqref{xhk} is Lotka-Volterra system, which has at most one center. And there exists a unique configuration of centers for Lotka-Volterra system in the equivalent sense, which is $(0;0;0;1)$. Therefore, the interesting problem is to study configurations of centers for system \eqref{xhk} with degree $n\ge 3$.
 We now state the main result in this section.
\begin{theorem}\label{thm-main2}
Suppose that system \eqref{xhk} has  $[\frac{n^2-1}{2}]$ centers with $n\geq 3$ and its configurations are $(i_1;i_2;i_3;i_4)$. Then the following statements hold.
\begin{itemize}
\item[(i)]$i_j\neq 0,\quad j=2,3,4$.
\item[(ii)]$\min\{i_1+i_3,i_2+i_4\}\geq n-1$ (resp. $n-2$) if $n$ is odd (resp. even).
\item[(iii)]If $i_1=0$, then $i_2\geq[\frac{n-1}{2}]$ and $i_3\geq n-2$ (resp. $n-1$) when $n$ is even (resp. odd).
\item[(iv)]$[\frac{n^2+4}{8}]\leq \max\{i_1,i_2,i_3,i_4\}\leq [\frac{n^2-2n+2}{2}]$ and $i_1\leq[\frac{n^2-1}{8}]$.
\item[(v)]If there exists $j\in\{1,2,3,4\}$ such that $i_j=[\frac{n^2-2n+2}{2}]$, then the configuration of centers of  system \eqref{xhk} must be
    $$
    (i_1;i_2;i_3;i_4)=\left(0;[\frac{n-1}{2}];[\frac{n^2-2n+2}{2}];[\frac{n-1}{2}]\right).
    $$
\end{itemize}
\end{theorem}

Before proving Theorem \ref{thm-main2}, let us give some preliminaries. Note that the set consisting of HK-vector fields $X_{hk}\in \Psi_{n,n}$ is generic in the space of all HK-vector fields $X_{hk}$. Using the similar arguments in proof
  of Proposition \ref{pro-center}, we can obtain %it is no hard to prove HK-vector fields which belong to $\Psi_{n,n}$ is generic in HK-vector fields which belong to $\mathcal{Y}_{n,n}$. Hence, we have
\begin{proposition}\label{pro-1}
Suppose the HK-vector field $X_{hk}^*\in\mathcal{Y}_{n,n}\setminus\Psi_{n,n}$ has $[\frac{n^2-1}{2}]$ centers. Then there exists a HK-vector field ${X}_{hk}\in\Psi_{n,n}$ such that $ {X}_{hk}$ has the same configuration of centers to that of $X_{hk}^*$.
\end{proposition}
For  convenience, we denote the vector fields $X_{hk}\in \Psi_{n,n}$ of system \eqref{xhk} having $[\frac{n^2-1}{2}]$ centers by
$$
\mathcal{HK}=\left \{X_{hk}: \ X_{hk}\in \Psi_{n,n}, \ {\rm system}\ \eqref{xhk}\  {\rm has}\; [\frac{n^2-1}{2}] \ {\rm centers} \right\}.
$$
By Proposition \ref{pro-1}, we only need to consider the configurations of centers of $X\in \mathcal{HK}$. By Lemma \ref{lem-inf-node} and Lemma \ref{lem-equiv}, we have the following proposition.
\begin{proposition}\label{hkm}
Assume that vector fields $X\in \mathcal{HK}$, that is $X\in \Psi_{n,n}$ and  $X$ has $[\frac{n^2-1}{2}]$ centers. Then
\begin{itemize}
\item[(a)] when $n$ is even,  $X$ has $\frac{n^2-2}{2}$ elementary centers and $\frac{n^2+2}{2}$ elementary saddles in $\mathbb{R}^2$, and six  infinite critical points which are elementary nodes.
\item[(b)] when $n$ is odd, $X$ has $\frac{n^2-1}{2}$ elementary centers and $\frac{n^2+1}{2}$ elementary saddles in $\mathbb{R}^2$, and four infinite critical points which are elementary nodes.
\end{itemize}
\end{proposition}
To study $X\in \mathcal{HK}$ in the complete plane $\mathbb{R}^2$ including its behavior near infinity, it is suffice to study $X$ in the  Poincar\'e disc. This disc
  is divided into four sectors $\mathcal{S}_j$ by invariant lines $x=0$ and $y=0$ of $X$, that is,
  $$\mathcal{S}_j={\small \{(x,y):\  x=r\cos\theta,\ y=r\sin\theta, \ (j-1)\frac{\pi}{2}\le\theta\le j\frac{\pi}{2}, \ 0\le r\le M,\ M\gg 1  \}, }$$
  where $j=1,\cdots,4.$ Note that all the critical points on the boundaries $\theta=(j-1)\frac{\pi}{2}, 0<r<M$ and $\theta=j\frac{\pi}{2}, 0<r<M$ of $\mathcal{S}_j$ are elementary saddles and all the infinite critical points on the boundary $r=M$ of $\mathcal{S}_j$ are elementary nodes. Therefore,  the sum of indices at the critical points  of $X$ in the interior of $\mathcal{S}_j$, denoted by $\sum_{int_j}i$,  can be characterized by  the indices at the critical points on boundary of $\mathcal{S}_j$ as follows.
  \begin{lemma}\label{lem-poin-hop}
Suppose that $X\in \mathcal{HK}$, and $X$ has $s_j(X)$ saddles and $n_j(X)$ nodes in the boundaries of $\mathcal{S}_j$  except the three vertexes of $\mathcal{S}_j$: $(0,0)$, {\small $(M\cos((j-1)\frac{\pi}{2}),M\sin((j-1)\frac{\pi}{2}))$ and $(M\cos(j\frac{\pi}{2}),M\sin(j\frac{\pi}{2}))$}, $j=1,\cdots, 4$. Then
$$
s_j(X)=2\sum_{int_j}i+n_j(X), \ j=1, \cdots, 4.
$$
\end{lemma}

\begin{proof}
Since $X\in \mathcal{HK}$, by Proposition \ref{hkm} we know that either $n_j(X)=0$ or $n_j(X)=1$ if $n$ is even,  and $n_j(X)=0$ if $n$ is odd for $j=1,\cdots,4$.
For simplicity,  we let
$$
U(x,y)=-F(x,y)-y\frac{\partial F}{\partial y}, \ W(x,y)=F(x,y)+x\frac{\partial F}{\partial x}.
$$
 Consider the corresponding system \eqref{xhk} of  $X$ in the first sector $\mathcal{S}_1$, by the transformation $(x,y)\mapsto (\sqrt{x},\sqrt{y})$ we have
\begin{equation}\label{eq11}
\begin{split}
\dot{x}&=\frac{1}{2}xU(x^2,y^2),\\
\dot{y}&=\frac{1}{2}yW(x^2,y^2).
\end{split}
\end{equation}
Since the systems \eqref{xhk} and \eqref{eq11} are topologically conjugate in the sector $\mathcal{S}_1$,  the number and topological classification of finite and infinite critical points of system \eqref{eq11} and system \eqref{xhk}  are the same in the sector $\mathcal{S}_1$.  Note that system \eqref{eq11} is invariant under transformations $(x,y)\mapsto(-x,y)$ $(x,y)\mapsto(x,-y)$ and $(x,y)\mapsto(-x,-y)$. Therefore,  system \eqref{eq11} can be defined in $\mathbb{R}^2$, which  has  $4n_1(X)+4$ nodes at infinity, and $2s_1(X)+1$ saddles on $x$-axis and $y$-axis. Applying Poincar\'{e}-Hopf theorem to system \eqref{eq11}, we obtain that
$$
2\sum_f i+\sum_{inf}i=2\left(-(2s_1(X)+1)+4\sum_{int_1}i\right)+4n_1(X)+4=2,
$$
and it follows $s_1(X)=2\sum_{int_1}i+n_1(X)$.

Using the similar arguments it can be  prove that $s_j(X)=2\sum_{int_j}i+n_j(X)$ for $j=2,3,4$. We omitted them to save space. The proof is finished.
\end{proof}
Lemma \ref{lem-poin-hop} is a powerful tool to study the configuration of centers for $X\in \mathcal{HK}$. It tell us that the configurations of centers can be controlled by the configurations of all the saddles of $X$ as follows.
\begin{corollary}\label{coro1}
Suppose that $X\in \mathcal{HK}$. If there exists a sector $\mathcal{S}_j$ such that  $X$ has at least two (resp. one) saddles except the origin $(0,0)$ on the $x$-axis and $y$-axis in $\mathcal{S}_j$ when $n$ is even (resp. odd), then there exists at least one center in the interior of this sector.
\end{corollary}
Now we are ready to prove Theorem \ref{thm-main2}.

\begin{proof}[Proof of Theorem \ref{thm-main2}] Due to Proposition \ref{pro-1}, we only need to prove Theorem \ref{thm-main2} for $X\in \mathcal{HK}$.
 It is clear that conclusion (i) of Theorem \ref{thm-main2} comes  directly from Proposition \ref{hkm} and Corollary \ref{coro1}.

 Note that the  arguments applied to verify conclusions (ii) - (v) are similar for both even $n$ and odd $n$. Thus, in the following we consider only the case that $n$ is even.

We first prove conclusion (ii). By Lemma \ref{lem-f}, we have
$$
i_{j}\geq \sum_{int_j}i,
$$
where $i_j$ is the number of centers in the $j$th sector.
Notes that $X\in \mathcal{HK}$, which has $2(n-1)$ saddles on the set $\{(x,y)\neq(0,0):xy=0\}$ and six nodes at infinity. By Lemma \ref{lem-poin-hop}, we know that one of following statements holds:
\begin{itemize}
\item[(i)]$\sum_{int_1}i+\sum_{int_3}i=\sum_{int_2}i+\sum_{int_4}i+1=n-1$, when $F_{n-1}(x,y)$ has a linear factor $x+ky$ with $k>0$;
\item[(ii)]$\sum_{int_1}i+\sum_{int_3}i+1=\sum_{int_2}i+\sum_{int_4}i=n-1$, when $F_{n-1}(x,y)$ has a linear factor $x+ky$ with $k<0$,
\end{itemize}
where $F_{n-1}(x,y)$ is the $n-1$-th homogenous part of polynomial $F(x,y)$. Hence, conclusion (ii) holds.

Let us now  prove conclusion (iii). If there is no centers in the first sector $\mathcal{S}_1$, then the number of saddles on $\{(x,0):x>0\}\cup\{(0,y):y>0\}$ is less than one by Corollary \ref{coro1}. From some calculations in two cases that $F_{n-1}(x,y)$ has a linear factor $x+ky$ with $k>0$ and $k<0$, respectively, one of following statements holds
\begin{itemize}
\item[(i)] if $k>0$, $i_1=\sum_{int_1}i=0,\ \sum_{int_3}i=n-1,\ \sum_{int_2}i=\sum_{int_4}i=\frac{n-2}{2}$.
\item[(ii)] if $k<0$, $i_1=\sum_{int_1}i=0,\ \sum_{int_3}i=n-2$ and either $\sum_{int_2}i=\frac{n}{2}$, $\sum_{int_4}i=\frac{n-2}{2}$ or $\sum_{int_2}i=\frac{n-2}{2}$, $\sum_{int_4}i=\frac{n}{2}$.
\end{itemize}
So we have
$$
i_2\geq\sum_{int_2} i\geq\frac{n-2}{2}=[\frac{n-1}{2}],\quad i_3\geq \sum_{int_3}i\geq n-2
$$
when $n$ is even. That is conclusion (iii).

To prove conclusion (iv), we assume that the number of centers of $X$ is the maximum in some  $\mathcal{S}_j$ of four sectors,  without loss of generality, we assume that
 $$i_3=\max\{i_1,i_2,i_3,i_4\}.$$
 Then
\begin{equation}\label{ineq}
\begin{split}
i_3 & \leq i_1+i_3=\frac{n^2-2}{2}-(i_2+i_4)\leq \frac{n^2-2}{2}-(\sum_{int_2}i+\sum_{int_4}i)\\
&\leq \frac{n^2-2}{2}-(n-2)=\frac{n^2-2n+2}{2}.
\end{split}
\end{equation}
On the other hand,  there exists a sector $\mathcal{S}_j$ such that the number $i_j$ of centers in $\mathcal{S}_j$ satisfying  $i_j\geq\frac{[\frac{n^2-1}{2}]}{4}$. Otherwise
$$[\frac{n^2-1}{2}]=\sum_{j=1}^4i_j< 4\cdot\frac{[\frac{n^2-1}{2}]}{4}=[\frac{n^2-1}{2}], $$ which is a contradiction. We now claim that
$$
\left\{\frac{[\frac{n^2-1}{2}]}{4}\right\}=[\frac{n^2+4}{8}],
$$
where $\{m\}$ represents {\it the minimum integer which is not less than $m$}.

In fact, since $n$ is even, there exists an integer $k$ such that either $n=4k$ or $n=4k+2$. We have
\begin{equation*}
\begin{split}
\left\{\frac{[\frac{(4k)^2-1}{2}]}{4}\right\}&=\left\{\frac{(4k)^2-2}{8}\right\}=2k^2=[\frac{(4k)^2+4}{8}], \ \text{as}\ n=4k.\\
\left\{\frac{[\frac{(4k+2)^2-1}{2}]}{4}\right\}&=\left\{\frac{16k^2+16k+2}{8}\right\}\\
&=2k^2+2k+1=\frac{(4k+2)^2+4}{8}, \ \text{as}\ n=4k+2.
\end{split}
\end{equation*}
Hence, $[\frac{n^2+4}{8}]\le i_3$.
Using the same method, we can obtain that $i_1\leq[\frac{n^2-1}{8}]$.

 Finally, we prove (v). By inequality \eqref{ineq}, $i_1=0$ if $i_3=\frac{n^2-2n+2}{2}$.

 Since $i_1=\min\{i_1,i_2,i_3,i_4\}$, we can always assume that $i_3=\frac{n^2-2n+2}{2}$. Note that inequality \eqref{ineq} are actually an equality. Thus, we have
 $$
 i_2+i_4=\sum_{int_2}i+\sum_{int_4}i=n-2,\quad i_1=0.
 $$
 Then we have following statements.
\begin{itemize}
\item[(a)]$\sum_{int_2}i+\sum_{int_4}i=n-2$, which  implies that $F_{n-1}(x,y)$ has a linear factor $x+ky$ with $k>0$.
\item[(b)]$i_2+i_4=\sum_{int_2}i+\sum_{int_4}i$, which  implies that there is no saddles in the interior of $\mathcal{S}_2$ and $\mathcal{S}_4$. Hence, $i_2=\sum_{int_2} i$ and $i_4=\sum_{int_4} i$.
\item[(c)] By Lemma \ref{lem-poin-hop} and statement (a), $i_1=0$  which implies all the saddles on $x$-axis and $y$-axis lie in $\{(x,0):x\leq 0\}\cup\{(0,y):y\leq 0)\}$. It can be calculated by Lemma \ref{lem-poin-hop} that
$$\sum_{int_1} i=0,\; \sum_{int_2} i=\frac{n-2}{2},\; \sum_{int_3} i=n-1,\; \sum_{int_4} i=\frac{n-2}{2}.$$
\end{itemize}
From (b) and (c), we have $i_2=i_4=\frac{n-2}{2}$. Thus,  $X$ has the configuration
    $$
    (i_1;i_2;i_3;i_4)=(0;\frac{n-2}{2};\frac{n^2-2n+2}{2};\frac{n-2}{2})=(0;[\frac{n-1}{2}];[\frac{n^2-2n+2}{2}];[\frac{n-1}{2}]).
    $$
Hence,  Theorem \ref{thm-main2} is verified.
\end{proof}

\section{Dynamics of cubic polynomial Kolmogorov systems with the maximum centers}\label{app}

In this section, we study the dynamics  of cubic  polynomial Kolmogorov vector fields $Y_k$ having four centers,  where $Y_{k}=(xP(x,y), yQ(x,y))$, $P(x,y)$ and $Q(x,y)$ are any two quadratic polynomials.
Especially, if the cubic polynomial Kolmogorov vector fields are Hamiltonian, authors in \cite{LX} have systematically investigated the configurations of centers for  the number and type of all possible centers, and left an open question
{\it  if there are only two types configurations of centers when the cubic polynomial Hamiltonian Kolmogorov vector fields $X_{hk}$ have four centers}.

 Applying  Theorem \ref{thm-main2}, we answer this open question affirmatively in the sense of equivalence, and obtain all global phase portraits for this cubic vector fields $X_{hk}$ having four centers.

 If the cubic polynomial Kolmogorov vector fields $Y_k$ are not Hamiltonian, we show that the cubic vector fields $Y_k$ have a first integral, which is well-defined elementary function on $\mathbb R^2$ except the set $\{(x,y): xy=0\}$, and there exist only three types of configurations of centers for the $Y_k$ in the equivalent sense. This reveals the difference between Hamiltonian and non-Hamiltonian integrable systems.

 \subsection{Cubic polynomial Hamiltonian Kolmogorov systems} 
\begin{figure}[htbp]
%\centering
%\subfigure{
%\begin{minipage}[t]{0.45\linewidth}
\centering
{\includegraphics[width=2in]{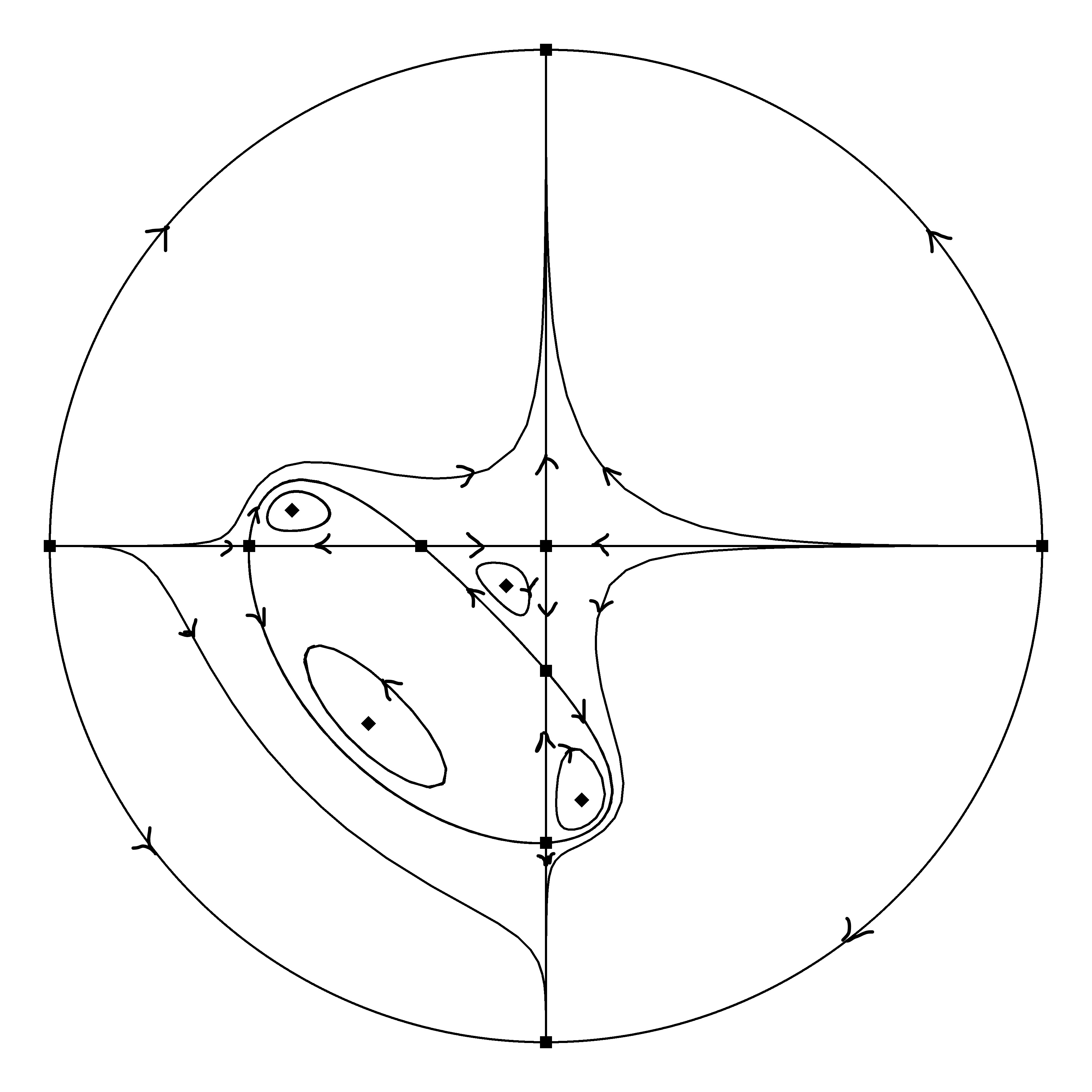}} %\setcounter{figure}{0}
\caption{\small Global phase portrait of system \eqref{3hk} with the configuration $(0;1;2;1)$ of centers}
\label{fig0121}
\end{figure}

\begin{figure}[htbp]
%\subfigure{
%\begin{minipage}[t]{0.45\linewidth}
\centering
{\includegraphics[width=2in]{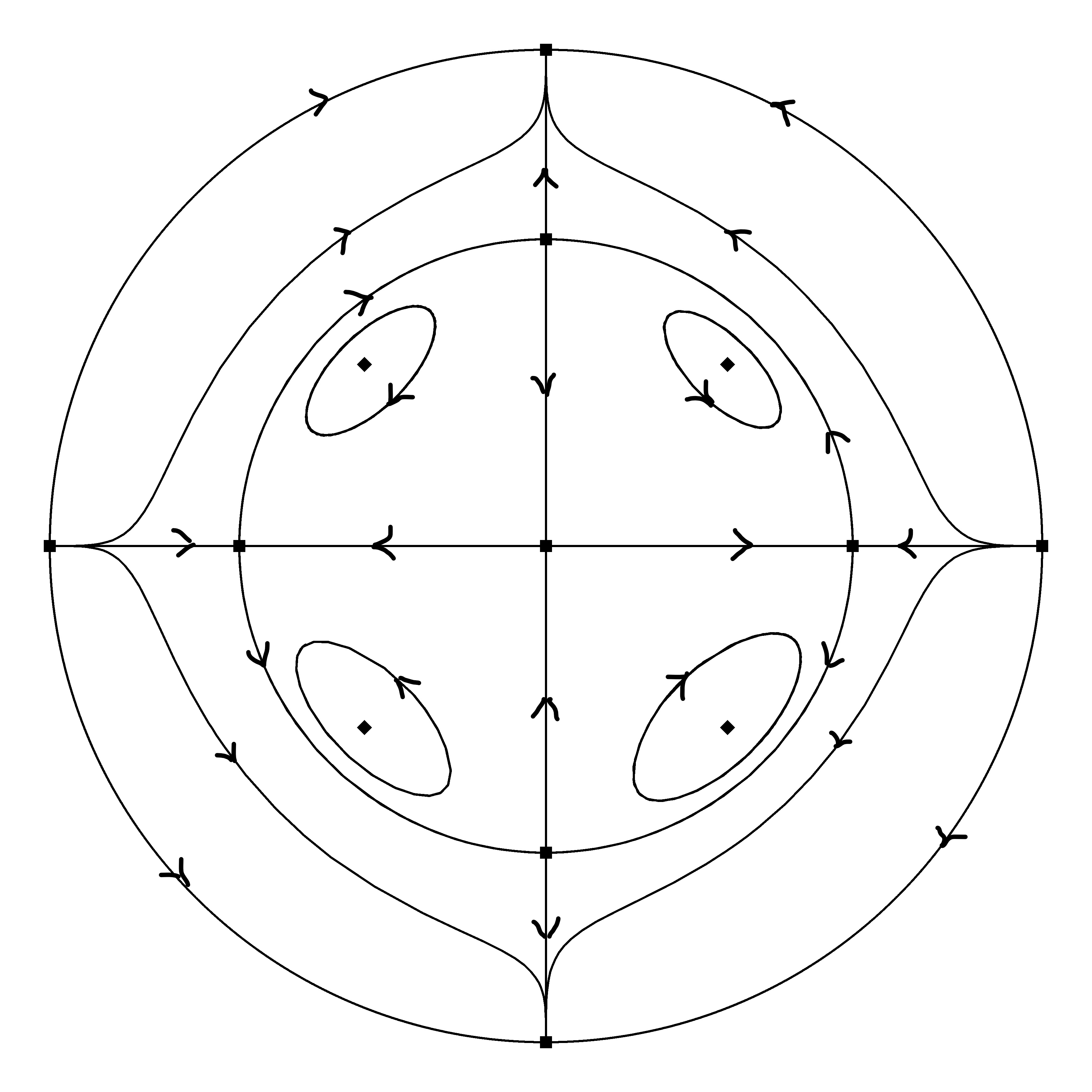}}
\caption{\small Global phase portrait of system \eqref{3hk} with the configuration $(1;1;1;1)$ of centers}
\label{fig1111}
%\end{minipage}}%}
\end{figure}

Consider the corresponding system of cubic Hamiltonian vector fields $X_{hk}$
 \begin{equation}\label{3hk}
\begin{split}
\frac{dx}{dt}=&-x\left(\hat F_2(x,y)+y\frac{\partial \hat F_2}{\partial y}\right),\\
\frac{dy}{dt}=& y\left(\hat F_2(x,y)+x\frac{\partial \hat F_2}{\partial x}\right),
\end{split}
\end{equation}
where $\hat F_2(x,y)$ is any a quadratic polynomial.  Clearly system \eqref{3hk} has at most four centers.
  Authors in \cite{LX} have founded two configurations $(1;1;1;1)$ and $(1;0;1;2)$ of centers if system \eqref{3hk} has four centers. According to the equivalence of two configurations in section 4, it can be checked that the configuration $(1;0;1;2)$ of centers is equivalent to  $(0;1;2;1)$. Thus, the open question proposed in \cite{LX} is ask {\it if there are only two configurations  $(1;1;1;1)$ and $(0;1;2;1)$ of centers in the sense of equivalence when system \eqref{3hk} has four centers}. In the following we
answer this open question affirmatively and obtain  the global phase portraits of system \eqref{3hk}.
\begin{theorem}\label{3hkT}
Suppose that system \eqref{3hk} has four centers. Then system \eqref{3hk} has  only two configurations $(0;1;2;1)$ and $(1;1;1;1)$ of centers in the sense of equivalence. Furthermore, the global dynamics of system \eqref{3hk} can be characterized as follows.
\begin{itemize}
\item[(i)]System \eqref{3hk} has exactly nine finite critical points, in which four centers and five saddles, and all five saddles lie in the level set $xy\hat F_2(x,y)=0$.
\item[(ii)]System \eqref{3hk} has exactly two pairs of infinite critical points which correspond to infinity in the direction of the $x$-axis and $y$-axis. Those infinite critical points are all elementary nodes.
\item[(iii)]Using linear transformations \eqref{trans} of  variables $(x,y)$ and time change $t\mapsto -t$ if necessary, system \eqref{3hk} has only two different topological types of global phase portraits in the Poincar\'e disc, which are sketched in Figure \ref{fig0121} and Figure \ref{fig1111}.
\end{itemize}
\end{theorem}

%\end{minipage}}
\begin{proof}
\setcounter{figure}{2}

Since $n=3$  and system \eqref{3hk} has four centers,   by conclusion (i) in Theorem \ref{thm-main2}, we know that system \eqref{3hk} has at least one center in $\mathcal{S}_j$ for $j=2,3,4$, that is, $i_2\ge 1$, $i_3\ge 1$ and $i_4\ge 1$.

If system \eqref{3hk} has one center in $\mathcal{S}_1$, that is $i_1=1$, then there exists a unique configuration $(1;1;1;1)$ of four centers for system \eqref{3hk}.

If system \eqref{3hk} has no centers in $\mathcal{S}_1$, that is $i_1=0$, then $i_3\ge 3-1=2$ and $i_2+i_4\ge 3-1=2$  by conclusion (ii) in Theorem \ref{thm-main2}. Thus, $i_2=1,i_3=2,i_4=1$. This implies system \eqref{3hk} has a unique configuration $(0;1;2;1)$
of four centers.

Summarizing the above analysis, we obtain that system \eqref{3hk} has only two configurations $(1;1;1;1)$ and $(0;1;2;1)$ of four centers.

We now discuss the global dynamics of system \eqref{3hk} with four centers.  Since system \eqref{3hk} has four centers, the following equations
$$
\hat F_2(x,y)+y\frac{\partial \hat F_2}{\partial y}=0, \ \hat F_2(x,y)+x\frac{\partial \hat F_2}{\partial x}=0
$$
have four solutions in the interior of  $\mathcal{S}_j$ for some $j\in\{1,2,3,4\}$. Thus,
 the quadratic  homogenous parts of polynomials
 $$\hat F_2(x,y)+y\frac{\partial \hat F_2}{\partial y}\ \ {\rm and}\ \ \hat F_2(x,y)+x\frac{\partial \hat F_2}{\partial x}$$
has no common linear factors. This implies that  $-x(\hat F_2(x,y)+y\frac{\partial \hat F_2}{\partial y})$ and $y(\hat F_2(x,y)+x\frac{\partial \hat F_2}{\partial x})$ has no common linear factors. Hence, the vector fields $X_{hk}$ of system \eqref{3hk} belongs to $\Psi_{3,3}$ and $\mathcal{HK}$.

By conclusion (b) in Proposition \ref{hkm}, we obtain that system \eqref{3hk} has exactly nine finite critical points, in which four centers and five saddles, and four infinite critical points which are all elementary nodes.

On the other hand, by Lemma \ref{lem-inf-node}, $\hat F_2(x,y)$ has no real linear factors. Otherwise,  $xy\hat F_2(x,y)$ has three different real linear factors, then system \eqref{3hk} has at most $[\frac{3^2+1-3}{2}]=3$ centers by Theorem \ref{thmmain1}. This contradicts to the
fact system \eqref{3hk} has four centers. Therefore,  $xy\hat F_2(x,y)$ has only two different real linear factors $x$ and $y$, which corresponds to four infinite critical points in the direction of the $x$-axis and $y$-axis. And  it can be seen that five saddles lie in the level set of $xy\hat F_2(x,y)=0$. It follows that the conclusions (i) and (ii) hold.

Note that $\hat F_2(x,y)$  has no real linear factor. This implies that the quadratic curve $\hat F_2(x,y)=0$ is an ellipse. Saddle $(0,0)$ is a common point of two lines $x=0$ and $y=0$, and the other four saddles are common points of the ellipse and two lines $x=0$ and $y=0$, respectively.
There are only two possibilities that an ellipse has four common points with two lines $x=0$ and $y=0$, which correspond to the configurations $(0;1;2;1)$ and $(1;1;1;1)$ in the sense of equivalence. In any case, there are four compact
 regions whose boundary are all formed by $x=0$, $y=0$ or $\hat F_2(x, y)=0$. On the boundary of any compact region, $xy\hat F_2(x, y)\equiv 0$.  Thus in the
  interior of any compact region, there is one extreme point, which must be a center. Except these four centers and  five saddles, there is no other critical points. Hence one can easily obtain that there are only two different topological types of global phase portraits of system \eqref{3hk}  in the Poincar\'e disc by using linear transformations \eqref{trans} of  variables $(x,y)$ and time change $t\mapsto -t$ if necessary, see Figure \ref{fig0121} and Figure \ref{fig1111} respectively.
\end{proof}

{\bf Remark 4.1} Even system \eqref{3hk} has four centers, it can be verified that the level set $xy\hat F_2(x,y)=h$, for any $h\in \mathbb{R}$, does not have four ovals in $\mathbb{R}^2$, and there exist some $h_0$, $0\not=|h_0|$ such that the curves $xy\hat F_2(x,y)=h_0$ has three ovals in $\mathbb{R}^2$ for Figure \ref{fig0121} and two ovals in $\mathbb{R}^2$ for Figure \ref{fig1111}, respectively.

%In fact, it can be checked that the algebraic curve $xy\hat F_2(x,y)-h_0=0$ which is an invariant for system \eqref{quarticsystem}.

 \subsection{Cubic polynomial Kolmogorov systems} Consider the corresponding system of cubic Kolmogorov vector fields $Y_{k}$
 \begin{equation}\label{3k}
\begin{split}
\frac{dx}{dt}=& xP(x,y), \\
\frac{dy}{dt}=& yQ(x,y),
\end{split}
\end{equation}
where $P(x,y)$  and $Q(x,y)$ are any two quadratic polynomials. Notice that $x=0$ and $y=0$ are two invariant lines of system \eqref{3k}. Then the centers of system \eqref{3k} should be in the interior of $\mathcal{S}_j$ for some $j=1, 2,3,4$ if system \eqref{3k} has centers.
A natural question is to ask whether system \eqref{3k} has limit cycles if the number of its centers is maximum four, and
how many configurations of centers system \eqref{3k} has if it is not Hamiltonian and it has four centers. We completely answer the two questions in the subsection as follows.
\begin{theorem}\label{Confcenter}
If system \eqref{3k} has four centers, then it has an elementary first integral and no limit cycles. All the possible
configurations of centers are $(1;1;1;1)$, $(0;1;2;1)$ and $(0;1;1;2)$ in the sense of equivalence.
\end{theorem}

%Before proving the theorem we discuss the integrability of the system　We first study
In order to prove Theorem \ref{Confcenter}, we first  study the integrability of system \eqref{3k} if it has four centers.

\begin{proposition}\label{integrability}
If system \eqref{3k} has four centers, then system \eqref{3k} is integrable, that is, there exists an elementary first integral of system \eqref{3k} in $\mathbb{R}^2\setminus\{xy=0\}$. So system \eqref{3k} has no limit cycles.
\end{proposition}
\begin{proof}
If system \eqref{3k} has four centers at $p_i=(x_i,y_i)$, $i=1,2,3,4$, then $x_iy_i\not=0$,
$$ P(x_i, y_i)=0, \ Q(x_i, y_i)=0, \ i=1,2,3,4.$$
Hence, the quadratic algebraic curves $P(x, y)=0$ and $Q(x, y)=0$ have and only four intersection points, whose multiplicity is one by B\'ezout's Theorem.
Note that  the divergence of system \eqref{3k} at the center $p_i=(x_i,y_i)$ is zero, that is,
\begin{equation}\label{div}
 P(x_i,y_i)+x_i\frac{\partial P(x_i,y_i)}{\partial x}+Q(x_i,y_i)+y_i\frac{\partial Q(x_i,y_i)}{\partial y}=0, \ 1\le i\le 4.
\end{equation}
Therefore, according to Max Noether Fundamental Theorem (see \cite{Ful} for detail), the quadratic polynomial $P(x,y)+x\frac{\partial P(x,y)}{\partial x}+Q(x,y)+y\frac{\partial Q(x,y)}{\partial y}$ can be linearly represented by polynomials $P(x,y)$ and $Q(x,y)$,
in other words, there exist real constants $\alpha$ and $\beta$ such that for $\forall (x,y)\in \mathbb{R}^2$,
$$P(x,y)+x\frac{\partial P(x,y)}{\partial x}+Q(x,y)+y\frac{\partial Q(x,y)}{\partial y}=(1-\alpha) P(x,y)+(1-\beta) Q(x,y).$$

It is easy to check that $x^{\alpha-1}y^{\beta-1}$ is an integrable factor of system \eqref{3k} in $\mathbb{R}^2$ except $x$-axis and $y$-axis, that is, there is a function $\mathcal{F}(x, y)$   such that
$$\frac{\partial \mathcal{F}(x,y)}{\partial y}=-x^\alpha y^{\beta-1}P(x, y), \quad \frac{\partial \mathcal{F}(x,y)}{\partial x}=x^{\alpha-1} y^{\beta}Q(x, y).$$
This function $\mathcal{F}(x, y)$ is called {\it a first integral} of system \eqref{3k}.

Suppose that
$$P(x, y)=\sum_{i+j=0}^2p_{ij}x^iy^j, \quad Q(x, y)=\sum_{i+j=0}^2q_{ij}x^iy^j, \ p_{ij},q_{ij}\in \mathbb{R}, \ i,j\in\{0,1,2\}.$$
Then
\begin{equation}\label{coePQ} (\alpha+i)p_{ij}+(\beta+j)q_{ij}=0, \quad 0\leq i+j\leq 2.\end{equation}
We now discuss the form of $\mathcal{F}(x,y)$ depending on the values of $\alpha$ and $\beta$.

If $\alpha\not\in \{0, -1, -2\}$, then $\alpha+i\not=0$ since $0\le i\le 2$. Thus, the first integral of  system \eqref{3k} is
\begin{equation}\label{IF1}
\mathcal{F}(x, y)=x^\alpha y^{\beta}R(x, y), \quad R(x, y)=\sum_{i+j=0}^2\frac{q_{ij}}{\alpha+i}x^iy^j.
\end{equation}

Similarly, if $\beta\not\in \{0, -1, -2\}$, then $\beta+j\not=0$ since $0\le j\le 2$. So the first integral of  system \eqref{3k} is
\begin{equation}\label{IF2}
\mathcal{F}(x, y)=x^\alpha y^{\beta}R(x, y), \quad R(x, y)=-\sum_{i+j=0}^2\frac{p_{ij}}{\beta+j}x^iy^j.
\end{equation}

If $\alpha+\beta\not\in \{0, -1, -2\}$, then $(\alpha+i)^2+(\beta+j)^2\not=0$. Hence, the first integral of  system \eqref{3k} is
\begin{equation}\label{IF3}
\begin{split}
\mathcal{F}(x, y)&=x^\alpha y^{\beta}R(x, y), \\
 R(x, y)&=\sum_{i+j=0, \alpha+i\not=0}^2\frac{q_{ij}}{\alpha+i}x^iy^j-\sum_{i+j=0, \alpha+i=0}^2\frac{p_{ij}}{\beta+j}x^iy^j.
 \end{split}
\end{equation}

If $\alpha, \beta, \alpha+\beta\in \{0, -1, -2\}$, then the first integral of  system \eqref{3k} may contain logarithmic functions.
Concretely, if $(\alpha, \beta)\in \{(0, 0), (0, -1), (0, -2), (-1, -1)\}$, then the first integrals of  system \eqref{3k} are
\begin{equation}\label{IF4}
\begin{array}{ll}
\mathcal{F}(x, y)=-p_{00}\ln |y|-(p_{01}+p_{11}x)y-\frac{p_{02}y^2}{2}+q_{00}\ln |x|+q_{10}x+\frac{q_{20}x^2}{2},\\
\mathcal{F}(x, y)=-p_{01}\ln |y|+y^{-1}(p_{00}+p_{10}x+p_{20}x^2)-p_{02}y+q_{11}x+q_{01}\ln |x|,\\
\mathcal{F}(x, y)=\frac{y^{-2}}{2}(p_{00}+p_{10}x+p_{20}x^2)+p_{11}xy^{-2}-p_{02}\ln |y|+q_{02}\ln |x|,\\
\mathcal{F}(x, y)=x^{-1}y^{-1}(p_{00}+p_{10}x+p_{20}x^2-p_{02}y^2)-p_{11}\ln |y|+q_{11}\ln |x|,\\
\end{array}
\end{equation}
in $\mathbb{R}^2\setminus (\{(x,y): x=0\}\cup \{(x,y): y=0\})$, respectively.

For the other $(\alpha, \beta)\in \{(-1, 0), (-2, 0)\}$, one can calculate the first integral of system \eqref{3k} easily by symmetry, we omit it here.
Hence, system \eqref{3k} is integrable with a first integral $\mathcal{F}(x,y)$ almost everywhere in $\mathbb{R}^2$.
\end{proof}
As we have already observed in Proposition \ref{integrability} that the first integral of system \eqref{3k} may contain polynomials, or rational functions or logarithmic functions.
The following lemma shows that it is not necessary to consider these system \eqref{3k}  with the first integral  containing logarithmic functions in study the configurations of centers.

\begin{lemma}\label{Log} Assume that system \eqref{3k} with four centers has a first integral  containing logarithmic functions, then there exists another system \eqref{3k}
with a first integral without logarithmic functions such that the two systems have the same configuration of centers.
\end{lemma}
\begin{proof}
Since their proof is similar, we only give the detailed proof in the case $(\alpha, \beta)=(0, 0)$. Furthermore, we assume
$p_{00}q_{00}\not=0$, else the problem becomes easier.
Hence, system \eqref{3k} with four centers has the form
\begin{equation}\label{00}
\begin{split}
\frac{dx}{dt}=& x(p_{00}+p_{01}y+p_{11}xy+p_{02}y^2), \\
\frac{dy}{dt}=& y(q_{00}+q_{10}x-p_{11}xy+q_{20}x^2),
\end{split}
\end{equation}
since system \eqref{3k} has the first integral $\mathcal{F}(x,y)$ with the first expression in \eqref{IF4} as $(\alpha, \beta)=(0, 0)$.

Consider a small perturbation of  system \eqref{00}
\begin{equation}\label{00p}
\begin{split}
\frac{dx}{dt}=& x(p_{00}+p_{01}y+p_{11}xy+p_{02}y^2)-\epsilon p_{00}x\left(q_{10}x-p_{01}y+\frac{q_{20}x^2}{2}-p_{11}xy-\frac{p_{02}y^2}{2}\right), \\
\frac{dy}{dt}=& y(q_{00}+q_{10}x-p_{11}xy+q_{20}x^2)+\epsilon q_{00}x\left(q_{10}x-p_{01}y+\frac{q_{20}x^2}{2}-p_{11}xy-\frac{p_{02}y^2}{2}\right),
\end{split}
\end{equation}
where $0<\epsilon\ll 1$.

The first integral of system \eqref{00p} is
\begin{equation}\label{fe}
\mathcal{F}(x, y)=\frac{y^{-\epsilon p_{00}}x^{\epsilon q_{00}}}{\epsilon}\left(1+\epsilon(q_{10}x-p_{01}y)+\epsilon (\frac{q_{20}x^2}{2}-p_{11}xy-\frac{p_{02}y^2}{2})\right).
\end{equation}
Since $\epsilon>0$ is sufficiently small, in a small neighborhood of  each center of system \eqref{00}, system \eqref{00p} has a critical point, which is a center or a focus.

Note that  the centers of system \eqref{00p} are in the  region $xy\not=0$ and in this region system \eqref{00p} has an analytic first integral $\mathcal{F}(x, y)$ in \eqref{fe}. Thus, this critical point of system \eqref{00p} must be a center. This implies
that system \eqref{00p} has also four centers,  and the configuration of centers of system \eqref{00p} is the same to that of system \eqref{00}. The proof is complete.
\end{proof}

From now on, we only consider that system \eqref{3k}  has a first integral with the form
\begin{equation}\label{HR}
\mathcal{F}(x, y)=x^\alpha y^\beta R(x, y), \quad R(x, y)=\sum_{i+j=0}^2r_{ij}x^iy^j.
\end{equation}
Define
$$
\begin{array}{ll} \aleph\triangleq &\{\mathcal{F}(x, y)\in \eqref{HR}:\  \mbox{each of three polynomials}\ r_{00}+r_{10}x+r_{20}x^2, \\
&r_{00}+r_{01}y+r_{02}y^2,
 r_{20}x^2+r_{11}xy+r_{02}y^2\ \mbox{does not have multiple}\\
 &\mbox{ factors and}\ \alpha\beta (\alpha+\beta+2) r_{00}r_{20}r_{02}\not=0.\}
\end{array}
$$

The next result shows that we only need to consider the first integral $\mathcal{F}(x, y)\in \aleph$ in study the configurations of centers of system \eqref{3k}.

\begin{lemma} If system \eqref{3k} has four centers, then there exists  another  system \eqref{3k} which has
 a first integral $\mathcal{F}(x,y)\in\aleph$ such that the two systems have the same configuration of centers.
\end{lemma}
\begin{proof}
From Lemma \ref{Log}, we only need to consider system \eqref{3k} with the first integral \eqref{HR}.

Firstly, we show $\alpha\beta(\alpha+\beta+2)\neq 0$ if system \eqref{3k} has four centers. In fact, from the first integral in \eqref{HR}, system \eqref{3k} has the form
\begin{equation}\label{HR2}
\begin{split}
\frac{dx}{dt}=&xP(x,y)=x(-\beta R(x,y)-y\frac{\partial R}{\partial y}(x,y)), \\
\frac{dy}{dt}=&yQ(x,y)=y(\alpha R(x,y)+x\frac{\partial R}{\partial x}(x,y)).
\end{split}
\end{equation}
If $\alpha\beta=0$, for example $\alpha=0$, then $Q(x,y)=x\frac{\partial R}{\partial x}(x,y)$. Thus, by B\'ezout Theorem, $P(x,y)=0, Q(x,y)=0$ has at most two isolated zeros in the interior of four quadrants in $\mathbb{R}^2$. That means system \eqref{3k} has at most two centers which is a contradiction with four centers.

If $\alpha+\beta+2=0$, we have
$$
Q_2(x,y)-P_2(x,y)=\alpha R_2(x,y)+\beta R_2(x,y)+x\frac{\partial R_2}{\partial x}+y\frac{\partial R_2}{\partial y}=(\alpha+\beta+2)R_2(x,y)=0,
$$
where $P_2(x,y)$,  $Q_2(x,y)$ and $R_2(x,y)$ are the quadratic homogenous part of polynomials $P(x,y)$,  $Q(x,y)$ and $R(x,y)$, respectively.
Therefore, $Q(x,y)-P(x,y)$ is a linear polynomial. Hence, $P(x,y)=Q(x,y)=0$ which is equivalent to $P(x,y)=Q(x,y)-P(x,y)=0$ has at most two isolated zeros. It is also a contradiction with the fact there are four centers.

We now consider the  case  $\mathcal{F}(x,y)\not\in \aleph$.  If $r_{00}r_{20}r_{02}=0$, then we use the similar arguments in the proof of Lemma \ref{Log} to find a system whose first integral is in $\aleph$.  For example, we consider  a small perturbation Kolmogorov system of system \eqref{3k} having the following first integral
$$
\mathcal{F}_\epsilon(x, y)=x^\alpha y^\beta (R(x, y)+\epsilon(1+x^2+y^2)),
$$
where $0<\epsilon\ll 1$.  This perturbation Kolmogorov system
 has the same configuration of centers to that of the original system and $(r_{00}+\epsilon)(r_{20}+\epsilon)(r_{02}+\epsilon)\neq0$.

Then we consider the system \eqref{3k} with $r_{00}r_{20}r_{02}\not=0$, but there is at least one of  $r_{00}+r_{10}x+r_{20}x^2, r_{00}+r_{01}y+r_{02}y^2$ or $r_{20}x^2+r_{11}xy+r_{02}y^2$ which has multiple factor, that is
$$
(r_{10}^2-4r_{00}r_{20})(r_{01}^2-4r_{00}r_{02})(r_{11}^2-4r_{20}r_{02})=0.
$$
Then using perturbation technical again, we consider a small perturbation of this system such that the perturbation system has the first integral
$$
\mathcal{F}_\epsilon(x, y)=x^\alpha y^\beta (R(x, y)+\epsilon(x+y+xy)),
$$
and has the same configuration of centers to that of the original system. For the perturbation Kolmogorov system, we have
$$
(r_{10}+\epsilon)^2-4r_{00}r_{20}\neq 0,\ (r_{01}+\epsilon)^2-4r_{00}r_{02}\neq 0,\ (r_{11}+\epsilon)^2-4r_{20}r_{02}\neq 0,
$$
This implies that for any system \eqref{3k} having four centers, there exists a Kolmogorov system such that this Kolmogorov system has a first integral $\mathcal{F}_\epsilon(x,y)\in\aleph$. The proof is complete.
\end{proof}

From now on, we only consider system \eqref{3k} with the first integral  $\mathcal{F}(x, y) \in \aleph$. We will determine the topological classification of the critical points of system \eqref{3k} on the $x$-axis and $y$-axis.

\begin{lemma}\label{lem-HR}
Suppose that system \eqref{3k} has a first integral $\mathcal{F}(x, y) \subset \aleph$. Then all the critical points and infinite critical points are elementary. Furthermore, the following conclusions hold.
\begin{itemize}
\item[(i)] If $\alpha<0$, then  all the critical points of system \eqref{3k} on the $y$-axis must be nodes except the origin $(0,0)$;
\item[(ii)] If $\beta<0$, then  all the critical points of system \eqref{3k} on the $x$-axis must be nodes except the origin $(0,0)$;
\item[(iii)] If $\alpha>0$ and $\beta>0$, then all the infinite critical points of system \eqref{3k} must be nodes;
\item[(iv)] If $\alpha\beta<0$, then the origin $(0,0)$ is a node; if $\alpha\beta>0$, then  the origin $(0,0)$ is a saddle.
\end{itemize}
\end{lemma}
\begin{proof} These conclusions can be proved by calculation of Jacobian matrix at corresponding critical points directly. Since the arguments are similar in proof of conclusions (i)-(iv), we only prove conclusion (i) to save the space.
We just mention the calculation on infinite critical points for conclusion (iii), which needs Poincar\'e compactification used in the proof of Theorem \ref{thmmain1}.

We now prove conclusion (i). If $\alpha<0$, assume that system \eqref{3k} has a critical point $p=(0, y^*)$ with $y^*\not=0$ on the $y$-axis, by equation \eqref{HR2}, then we have
$$
Q(0, y^*)=\alpha R(0, y^*)=0.
$$
Further, the Jacobian matrix of system \eqref{HR2} at point $p=(0,y^*)$ is
\begin{equation*}
\begin{bmatrix}
 -y^*\frac{\partial R}{\partial y}(0,y^*) & 0 \\
 (\alpha+1)y^*\frac{\partial R}{\partial x}(0,y^*) &  \alpha y^*\frac{\partial R}{\partial y}(0,y^*)
\end{bmatrix}.
\end{equation*}
Since $R(0,y)=r_{00}+r_{01}y+r_{02}y^2$ has no multiple factors, $y=y^*$ is a simple root of $R(0,y)$, i.e. $\frac{\partial R}{\partial y}(0,y^*)\neq 0$. It follows that $(0, y^*)$ is an elementary node.
\end{proof}

Using the same notations in section 4, we say system \eqref{3k} has the configuration $(i_1; i_2; i_3; i_4)$ of centers, which implies there are $i_j$ centers in the interior of $\mathcal{S}_j$, $i_j\ge 0$ and
 $i_1+i_2+i_3+i_4=4$. Since the first integral of system \eqref{3k} $\mathcal{F}(x, y) \in \aleph$, polynomials $R(x,0)=r_{00}+r_{10}x+r_{20}x^2,\  R(0,y)=r_{00}+r_{01}y+r_{02}y^2$ and $R_2(1,y)=r_{20}+r_{11}y+r_{02}y^2$ have $r_x^+, r_{y}^+, r_{inf}^+$ positive real roots and
 $r_x^-, r_{y}^-, r_{inf}^-$ negative real roots,  respectively. That means system \eqref{3k} has exactly $r_x^+$ critical points on positive $x$-axis and $r_x^-$ critical points on negative $x$-axis, $r_{y}^+$ critical points on positive $y$-axis and $r_y^-$ critical points on negative $y$-axis, $r_{inf}^++2$ infinite critical points in sector $\mathcal{S}_1$ and $r_{inf}^-+2$ infinite critical points in sector $\mathcal{S}_2$. Obviously, $r_x^++r_x^-, r_{y}^++r_{y}^-,r_{inf}^++r_{inf}^-\in \{0, 2\}$.

\begin{lemma} $\max\{i_1,i_2,i_3,i_4\}\leq 2$.
\end{lemma}

\begin{proof}
Note that  $\max\{i_1,i_2,i_3,i_4\}\in\{i_3,i_4\}$. Without loss of generality, let $i_3=\max\{i_1,i_2,i_3,i_4\}$. Consider the following  system
\begin{equation}\label{QCKS}
\begin{split}
\frac{du}{dt}&=-\frac{1}{2}uP(-u^2, -v^2), \\
\quad \frac{dv}{dt}&=-\frac{1}{2}vQ(-u^2, -v^2).
\end{split}
\end{equation}

System \eqref{QCKS} is topological conjugated with system \eqref{3k} in the sector $\mathcal{S}_3$ by the transformation $x=-u^2, y=-v^2$. Hence, all the finite
and infinite critical points of system \eqref{QCKS} and system \eqref{3k} have the same topological classification in $\mathcal{S}_3$.

On the other hand,  system \eqref{QCKS} is invariant under the transformation $(u, v)\mapsto (-u, v)$ or $(u, v)\mapsto (u, -v)$.
This implies that except the origin $(0,0)$, \eqref{QCKS} has $2r_x^-$ critical points on $x$-axis, $2r_{y}^-$ critical points on $y$-axis and
$4+4r_{inf}^+$ infinite critical points in $\mathbb{R}^2$.

If $\alpha<0, \beta<0$, then the origin is a saddle and the other critical points on $x$-axis and $y$-axis are all nodes. By Poincar\'{e}-Hopf theorem and Lemma \ref{lem-HR}, we have
$$
4i_3+2r_x^-+2r_y^--1=\sum_f i=\frac{1}{2}(2-\sum_{inf}i)\leq \frac{1}{2}(2+4+4r_{inf}^+)\leq 7,
$$
which implies that $i_3\leq 2$.

If $\alpha<0, \beta>0$, then the origin is a node and the critical points on $y$-axis are all nodes, thus
$$
4i_3-2r_x^-+2r_y^-+1\leq\sum_f i=\frac{1}{2}(2-\sum_{inf}i)\leq 7,
$$
which implies that $i_3\leq 2$. The case $\alpha>0, \beta<0$ follows from the symmetry.

If $\alpha>0$ and $\beta>0$, then  the origin is a saddle and all the infinite critical points are nodes, thus
$$
4i_3-2r_x^--2r_y^--1\leq\sum_f i=\frac{1}{2}(2-\sum_{inf}i)=\frac{1}{2}(2-(4+4r_{inf}^+))\leq -1,
$$
which implies that $i_3\leq 2$.
\end{proof}

\begin{lemma}\label{i1i2}
$i_3+i_4\leq 3$.
\end{lemma}

\begin{proof} Consider the system
\begin{equation}\label{QCKS2}
\begin{split}
\frac{du}{dt}&=\frac{1}{2}uP(u^2, -v^2), \\
\frac{dv}{dt}&=\frac{1}{2}vQ(u^2, -v^2).
\end{split}
\end{equation}

System \eqref{QCKS2} is topological conjugated with system \eqref{3k} in the sector $\mathcal{S}_4$ by the transformation $x=u^2, y=-v^2$.

Note that except the origin, system \eqref{QCKS2} has $2r_x^+$ critical points on $x$-axis, $2r_y^-$ critical points on $y$-axis and $4+4r_{inf}^-$ infinite critical points.

If $\alpha<0$ and $\beta<0$, by Poincar\'{e}-Hopf theorem and Lemma \ref{lem-HR},
$$4i_4-1\leq 4i_4+2r_x^++2r_y^--1=\sum_f i=\frac{1}{2}(2-\sum_{inf}i) \leq 3+2r_{inf}^-.$$
So $i_4\leq 1+\frac{1}{2}r_{inf}^-$. Similarly, by the same discussion on system \eqref{QCKS}, we have $i_3\leq 1+\frac{1}{2}r_{inf}^+$. Hence, we have
$$
i_3+i_4\leq 2+\frac{1}{2}(r_{inf}^++r_{inf}^-)\leq 3.
$$

If $\alpha<0$ and $\beta>0$, we have
$$4i_4-2r_x^++1\leq 4i_4-2r_x^++2r_y^-+1\leq \sum_f i=\frac{1}{2}(2-\sum_{inf}i) \leq 3+2r_{inf}^-.$$
Thus, $i_4\leq \frac{1}{2}+\frac{1}{2}r_{inf}^-+\frac{1}{2}r_x^+$. Similarly, we have $i_3\leq \frac{1}{2}+\frac{1}{2}r_{inf}^++\frac{1}{2}r_x^-$ from the discussion on system \eqref{QCKS}. Hence, we have
$$
i_3+i_4\leq 1+\frac{1}{2}(r_{inf}^++r_{inf}^-)+\frac{1}{2}(r_x^++r_x^-)\leq 3.
$$
From the symmetry, we can obtain that $i_3+i_4\leq 3$ in the case $\alpha>0, \beta<0$.

Last we consider the case: $\alpha>0$ and $\beta>0$. Then  we have
$$
4i_4-2r_x^+-5\leq4i_4-2r_x^+-2r_y^--1\leq\sum_f i=\frac{1}{2}(2-\sum_{inf}i)=-1-2r_{inf}^-\leq -1.
$$
Hence, $i_4\leq 1+\frac{1}{2}r_x^+$. Similarly,  we can obtain that  $i_3\leq 1+\frac{1}{2}r_x^-$. Thereby
$$
i_3+i_4\leq 2+\frac{1}{2}(r_x^++r_x^-)\leq 3.
$$
\end{proof}
Similar to proof of Lemma \ref{i1i2}, we have
\begin{lemma}\label{ili3}
 $\max\{i_1+i_3,i_2+i_4\}\leq 3$.
\end{lemma}

At last, we are in the position to prove Theorem \ref{Confcenter}.
\begin{proof}[Proof of Theorem \ref{Confcenter}]
Recall that  $i_1=\min\{i_1,i_2,i_3,i_4\}$,   $i_2\leq i_4$ and $i_1+i_2+i_3+i_4=4$. If $i_1=1$, then $i_1=i_2=i_3=i_4=1$. Therefore, system \eqref{3k} has the configuration $(1; 1; 1; 1)$ of centers.
In \cite{LX}, authors have given a HK system \eqref{3hk} which has the configuration $(1; 1; 1; 1)$ of centers. Here we give the following non-Hamiltonian Kolmogorov system
\begin{equation}\label{ex1}
\begin{split}
\frac{dx}{dt}&=x(1-x^2-3y^2), \\
\frac{dy}{dt}&=2y(-1+2x^2+y^2),
\end{split}
\end{equation}
which has the first integral
$$
\mathcal{F}_1(x,y)=x^2 y ( x^2 + y^2 -1),
$$
and the four centers $p_i=(x_i,y_i)$ of system \eqref{ex1}  are
$$(x_i, y_i)\thickapprox (0.63,0.45), (-0.63,0.45),\ (-0.63,-0.45),\ (0.63,-0.45), \ i=1,2,3,4.$$
See the first diagram in Figure \ref{kss}.

\begin{figure}
  \centering
  \includegraphics[width=1\textwidth]{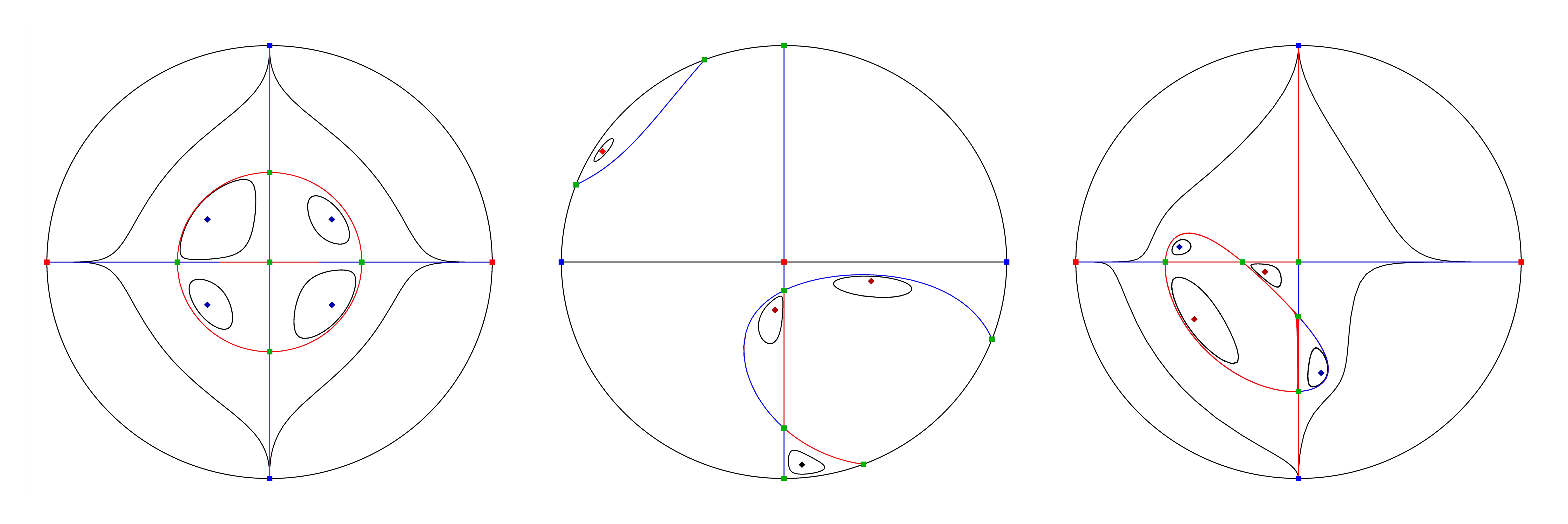}
  \caption{\small The configurations of centers by systems \eqref{ex1}, \eqref{ex2} and \eqref{ex3}, respectively.}\label{kss}
  \end{figure}

If $i_1=0$, then by Lemma \ref{i1i2} and Lemma \ref{ili3}, all the possible different configures of centers are $(0; 1; 1; 2)$ and $(0; 1; 2; 1)$.
Note that the configuration $(0; 1; 2; 1)$ of centers can be realized by HK system \eqref{3hk} in \cite{LX}, and Theorem \ref{3hkT} shows that the configuration  $(0; 1; 1; 2)$ of centers  can not be realized by HK system \eqref{3hk}.
However, we can find a non-Hamiltonian Kolmogorov system
\begin{equation}\label{ex2}
\begin{split}
\frac{dx}{dt}&=x(3-3x+8y+3x^2+6xy+y^2), \\
\frac{dy}{dt}&=\frac{1}{2}y(1-3x+4y+5x^2+9xy+y^2),
\end{split}
\end{equation}
which has the first integral
$$
\mathcal{F}_2(x,y)=\frac {({x}^{2}+3xy+{y}^{2}-x+4y+1) \sqrt {|x|}}{{y}^{3}}, \ \ \forall y\not=0,
$$
and the four centers $p_i=(x_i,y_i)$ of system \eqref{ex2}  are
$$(x_i, y_i)\thickapprox (-22.02,13.81),\  (-0.09, -0.47), \ (1.37,-15.94),\ (0.94, -0.21),\ i=1,2,3,4.$$
This implies that system \eqref{ex2} has the configuration  $(0; 1; 1; 2)$ of centers, see
 the second diagram in  Figure \ref{kss}. Thus, non-Hamiltonian Kolmogorov system \eqref{3k} has the configuration  $(0; 1; 1; 2)$ of centers.

Last we give an example to show that non-Hamiltonian Kolmogorov system \eqref{3k} has the configuration  $(0; 1; 2; 1)$ of centers as follows.
\begin{equation}\label{ex3}
\begin{split}
\frac{dx}{dt}&=-x(5 + 12 x+ 24 y+5 x^2 + 12 x y+ 15 y^2), \\
\frac{dy}{dt}&=y(15 + 48 x + 36 y + 25x^2 + 24 x y + 15 y^2),
\end{split}
\end{equation}
which has the first integral
$$
\mathcal{F}_3(x,y)=x^3 y (5 x^2 + 6 x y + 5 y^2 + 12 x + 12 y + 5),
$$
and the four centers $p_i=(x_i,y_i)$, $i=1,2,3,4$ of system \eqref{ex3}  are
$$(x_i, y_i)\thickapprox (-1.51,0.20),\ (-0.31,-0.09),\ (-1.31,-0.74),\ (0.28,-1.41).$$
See the third diagram in  Figure \ref{kss}. We finish the proof.
\end{proof}

{\bf Remark 5.2}: Theorem \ref{Confcenter} reveals the difference between cubic Hamiltonian Kolmogorov system \eqref{3hk} and cubic Kolmogorov system \eqref{3k}. For Hamiltonian Kolmogorov system \eqref{3hk}, we can obtain all different  global topological
phase portraits except the time reversal, see Figure \ref{fig0121} and \ref{fig1111}. However, we only give three  configurations of centers for cubic Kolmogorov system \eqref{3k}, see Figure \ref{kss}. It is interesting problem that how many different global phase portraits  the cubic Kolmogorov system \eqref{3k}
has, which is left for future study.

\section*{Acknowledgments}
Hongjin He and Dongmei Xiao are partially supported by National Key R $\&$ D Program of China (No. 2022YFA1005900), the Innovation Program of
Shanghai Municipal Education Commission (No. 2021-01-07-00-02-E00087) and the National Natural
Science Foundations of China (Nos. 11931016; 12271353). Changjian Liu is partially supported by the National Natural
Science Foundations of China (No. 12171491).

\end{document}